\documentclass[11pt]{article}
\usepackage{amssymb,amsmath,mathrsfs,amsfonts,graphicx,epsf} 
\usepackage{color}
\usepackage[dvipsnames]{xcolor}

 \input xy
 \xyoption{all}

 \usepackage{esint}

 \usepackage{fullpage}

\usepackage{authblk}

  \newcommand{\blue}{\color{blue}}

  \newcommand{\ffoot}[1]{}

\newcommand{\bD}{\blacktriangle}

\newcommand{\f}{\mathfrak{f}}

\usepackage{mathrsfs}
\usepackage{graphicx}
\usepackage{amssymb}
\usepackage{color}
\newtheorem{theorem}{Theorem}
\newtheorem{corollary}{Corollary}
\newtheorem{lemma}{Lemma}
\newtheorem{proposition}{Proposition}
\newtheorem{definition}{Definition}
\newtheorem{remark}{Remark}

\newcommand{\bt}{\begin{theorem}}
\newcommand{\et}{\end{theorem}}
\newcommand{\bl}{\begin{lemma}}
\newcommand{\el}{\end{lemma}}
\newcommand{\bp}{\begin{proposition}}
\newcommand{\ep}{\end{proposition}}
\newcommand{\bc}{\begin{corollary}}
\newcommand{\ec}{\end{corollary}}
\newcommand{\bdeff}{\begin{definition}}
\newcommand{\edeff}{\end{definition}}
\newcommand{\brem}{\begin{remark}}
\newcommand{\erem}{\end{remark}}

\newcommand{\con}{{\mathcal C}}

\newcommand{\bi}{\begin{itemize}}
\newcommand{\iii}{\item}
\newcommand{\ei}{\end{itemize}}
\newcommand{\bd}{\begin{description}}
\newcommand{\ed}{\end{description}}

\newcommand{\bqn}{\begin{eqnarray}}
\newcommand{\eqn}{\end{eqnarray}}
\newcommand{\eqnn}{\nonumber\end{eqnarray}}

\newcommand{\nn}{\nonumber}
\newcommand{\ba}[1]{\begin{array}{#1}}
\newcommand{\ea}{\end{array}}

\newcommand{\R}{\mathbb{R}}

\newcommand{\g}{\gamma}
\newcommand{\al}{\alpha}

\newcommand{\VecM}{\mathrm{Vec}(M)}

\newcommand{\Gq}{{\gg}_q}

\newcommand{\la}{\left\langle}
\newcommand{\ra}{\right\rangle}

\newcommand{\distr}{{\blacktriangle}}

\newcommand{\Zz}{\mathcal{Z}}
\renewcommand{\gg}{{\bf G}}

\newcommand{\Tone}{type-1}

\newcommand{\Ttwo}{type-2}
\xdefinecolor{og}{named}{OliveGreen}

\newcommand{\mc}[1]{\mathcal{ #1 }}
\newcommand{\EXP}{\mc{E}}
\newcommand{\tcut}{t_{cut}}
\DeclareMathOperator{\Cut}{Cut}
\newcommand{\tcon}{t_{con}}
\DeclareMathOperator{\Con}{Con}
\newcommand{\til}[1]{\widetilde{#1}}

\newcommand{\fact}[1]{#1\mathpunct{}!}

\newcommand{\thep}{\theta_\sigma^+}
\newcommand{\them}{\theta_\sigma^-}

\title{\LARGE \bf Local properties of
almost-Riemannian structures in dimension 3\thanks{This research has been supported  by the European Research Council, ERC
StG 2009 ``GeCoMethods", contract number 239748 and by the Fondation Math\'ematique Jacques Hadamard.}
}

\author[$\spadesuit$]{U.~Boscain} 
\affil[$\spadesuit$]{CNRS, CMAP \'Ecole Polytechnique, Palaiseau, France and
Team GECO, INRIA Saclay --
\^Ile-de-France

{\tt ugo.boscain@cmap.polytechnique.fr}}

\author[$\clubsuit$]{G.~Charlot}
\affil[$\clubsuit$]{Univ. Grenoble Alpes, IF, F-38000 Grenoble, 
CNRS, IF, F-38000 Grenoble, France
{\tt Gregoire.Charlot@ujf-grenoble.fr}}

\author[$\blacklozenge$]{M.~Gaye} 
\affil[$\blacklozenge$]{CMAP \'Ecole Polytechnique, Palaiseau, France, 
Team GECO, INRIA Saclay --
\^Ile-de-France 

{\tt moussa.gaye@cmap.polytechnique.fr}}

\author[$\heartsuit$]{P.~Mason} 
\affil[$\heartsuit$]{CNRS-LSS-Supélec, rue Joliot-Curie, 91192 Gif-sur-Yvette, France

{\tt Paolo.Mason@lss.supelec.fr}}

\begin{document}

\maketitle


\begin{abstract}

A 3D almost-Riemannian manifold is a generalized Riemannian manifold defined locally by 3 vector fields that play the role of an orthonormal frame, but could become collinear on some set $\Zz$ called the singular set. Under the Hormander condition, a 3D almost-Riemannian structure still has a metric space structure, whose topology is compatible with the original topology of the manifold. Almost-Riemannian manifolds were deeply studied in dimension 2.

In this paper we start the study of the 3D case which appear to be reacher with respect to the 2D case, due to the presence of abnormal extremals which define a field of directions on the singular set. We study the type of singularities of the metric that could appear generically, we construct local normal forms and we study abnormal extremals.
We then study the nilpotent approximation and the structure of the corresponding  small spheres.

We finally give some preliminary results about heat diffusion on such  manifolds.

\end{abstract}

\section{Introduction}

A $n$-dimensional Almost Riemannian  Structure ($n$-ARS for short) is a rank-varying sub-Riemannian structure  that can be  locally defined by a set of $n$ (and not by less than $n$) smooth vector fields  on a $n$-dimensional manifold, satisfying the H\"ormander condition  (see for instance \cite{AgrBarBoscbook,bellaiche,jean1,jean2}).   These vector fields play the role of an orthonormal frame. 

Let us denote by $\bD(q)$  the linear span of the  vector fields at a  point $q$.  
Around a point $q$ where  $\bD(q)$ is $n$-dimensional, the corresponding metric is Riemannian. 

On  the singular set  $\Zz$, where $\bD(q)$ has dimension dim$(\bD(q))< n$, the corresponding Riemannian metric is not well-defined. However, thanks to the H\"ormander condition, one can still define the Carnot--Caratheodory distance between two points, which happens to be finite and continuous.

Almost-Riemannian structures were deeply studied for $n=2$: they were introduced in the context of hypoelliptic operators \cite{baouendi,FL1,grushin1}; they appeared in  problems of  population transfer in quantum systems  \cite{q4,BCha,q1} and have applications to  orbital transfer in space mechanics  \cite{BCa,tannaka}.

For 2-ARS, generically, the singular set  is a $1$-dimensional embedded submanifold  (see \cite{ABS}) and there are three types of points: Riemannian points, Grushin points where   $\bD(q)$  is $1$-dimensional and 
$\dim(\bD(q)+[\bD,\bD](q))=2$ and tangency points where $\dim(\bD(q)+[\bD,\bD](q))=1$ and the missing direction is obtained with one more bracket. Generically tangency points are isolated (see \cite{ABS,crest}).

$2$-ARSs present very interesting phenomena. For instance, the presence of a singular set permits the conjugate locus  to be nonempty even if the Gaussian curvature is negative where it is defined (see \cite{ABS}). Moreover, a Gauss--Bonnet-type formula can be obtained \cite{ABS,high-order,euler}.

 In \cite{BCGS} a necessary and sufficient condition for two  2-ARSs  on the same compact manifold $M$ to be   Lipschitz equivalent was given.  In \cite{camillo} the heat and the Schr\"odinger equation with the Laplace--Beltrami operator on a 2-ARS were studied. In that paper it was proven that the singular set acts as a barrier for the heat flow and for a quantum particle, even though geodesics can pass through the singular set without singularities.
 
In this paper  we start the study of 3-ARSs. An interesting feature of these structures is that abnormal extremals could be present. 
Abnormal extremals are special candidates to be length minimisers that cannot be obtained as a solution of  Hamiltonian equations. They do not exist in Riemannian geometry and they could be present in sub-Riemannian geometry. Abnormal minimisers are responsible for the non sub-analyticity of the spheres in certain  analytic cases \cite{agrachev97} and they are the subject of one of the most important open question in sub-Riemannian geometry, namely ``are length minimisers  always smooth?'' (see  \cite{agra-open,monti-solo}).

Moreover the presence of abnormal minimisers seems related to the non analytic hypoellipticity of the sub-Laplacian (built as the square of the vector fields defining the structure  \cite{christ91,christ92,treves}).

The simplest example of analytic sub-Riemannian structure for which there are abnormal minimisers and for which the ``sum of the square'' is not analytic hypoelliptic is provided by a 3-ARS, namely the  Baouendi-Goulaouic example, defined by the following three vector fields:
\bqn
X_1(x,y,z)=\left(
\ba{c} 1\\0\\0\ea
\right),~~~
X_2(x,y,z)=\left(
\ba{c} 0\\1\\0\ea
\right),~~~
X_3(x,y,z)=\left(
\ba{c} 0\\0\\x\ea
\right).
\eqn
For this structure, the  trajectory $(0,y_0+t,0)$ is an abnormal minimizer and the Green functions of the operator 
\bqn
\partial_t\phi=\Delta\phi\mbox{ where }\Delta=X_1^2+X_2^2+X_3^2=\partial_x^2+ \partial_y^2 + x^2\partial_z^2
\eqn
are not analytic.

\begin{figure}
\begin{center}
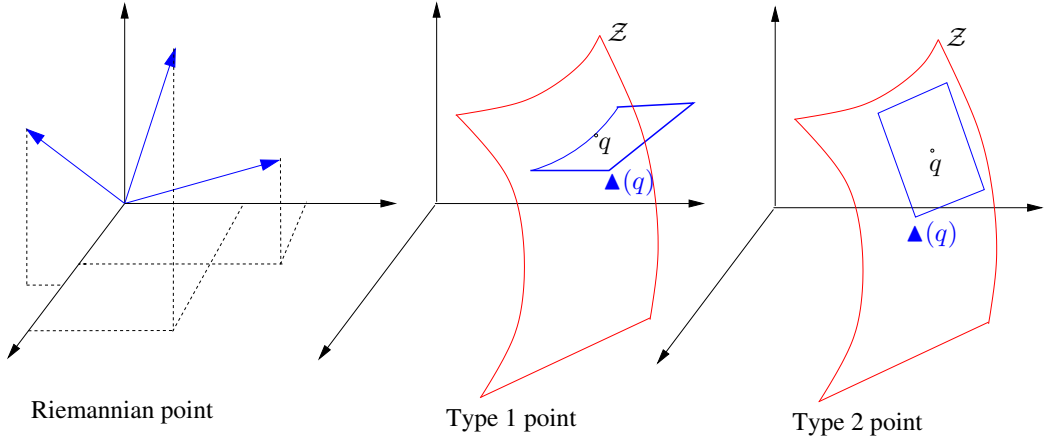
\caption{The 3 type of points occurring in the generic case\label{f-3-punti}}
\end{center}
\end{figure}

Our first result concerns the generic structure of singular sets for 3-ARSs. More precisely we prove that the following properties hold under generic conditions\footnote{ for the precise definition of generic see Definition \ref{d-generic}}
\bi
\item[(G1)] the dimension of $\bD(q)$ is larger than or equal to 2 and  $\bD(q)+[\bD,\bD](q)=T_qM$,  for every $q\in M$;
\item[(G2)] the singular set $\Zz$ (i.e., the set of points where dim$(\Zz)<3$) is an embedded 2-dimensional manifold;
\item[(G3)] the points where $\bD(q)=T_q\Zz$ are isolated.
\ei
As a consequence, under generic conditions, we can single out three types of points:
Riemannian points (where dim$(\bD(q))=3$), \Tone\ points  (where dim$(\bD(q))=2$ and $\bD(q)\neq T_q\Zz$)  and \Ttwo\ points (where $\bD(q)= T_q\Zz$). See Figure \ref{f-3-punti}. 
Moreover $\Zz$ is formed by \Tone\ and \Ttwo\ points (\Ttwo\ points are isolated).

Then for each type of point we build a local representation. See Theorem \ref{t-prenormal} below. These local representations are fundamental to study the local properties of the distance.

Next we study abnormal extremals. The intersection of $\bD$ with the tangent space to $\Zz$ define a field of directions on $\Zz$ whose integral trajectories are abnormal extremals. The singular points of this field of directions are \Ttwo\ points.

Since $\bD(q)+[\bD,\bD](q)=T_qM$ for every $q\in M$, as a consequence of  a theorem of A.Agrachev and J.P. Gauthier  \cite{AgrBarBoscbook,agra-gau-premio} if an abnormal extremal is optimal then it is not strict (i.e., it is at the same time a normal extremal). This condition is quite restrictive and indeed under generic conditions there are no abnormal minimisers.

We then focus on the nilpotent approximation of the structure at the different types of singular points. 
Obviously, at Riemannian points the nilpotent approximation is the Euclidean space.
Interestingly, at \Tone\ points the nilpotent approximation depends on the point and it is a combination of an Heisenberg sub-Riemannian structure and a Baouendi-Goulaouic 3-ARS. 
Surprisingly the nilpotent approximation at a \Ttwo\ point is the Heisenberg sub-Riemannian structure and hence is not a 3-ARS.

We then describe the metric spheres  for the nilpotent approximation for \Tone\ points.  As explained above, \Tone\ points are the only ones for which the nilpotent approximation is not yet well understood.
We recall that  for small radii, sub-Riemannian balls tend to those of the nilpotent approximation (after a suitable rescaling) in the Hausdorff-Gromov sense. 

We then study heat diffusion on 3-ARSs for nilpotent structures of \Tone. We study the diffusion related to two types of operators:
$\Delta_L$, the Laplace operator built as the divergence of the horizontal gradient, where the divergence is computed with respect to the Euclidean volume and $\Delta_\omega$, the Laplace operator built in the same way, but with  the divergence computed with respect to the intrinsic (diverging on $\Zz$) Riemannian volume. The first operator is hypoelliptic and for it we compute the explicit expression of the heat kernel. For the second operator we prove that it is not hypoelliptic and that the singular set $\Zz$ acts as a barrier for the heat flow.

\medskip
{\bf Structure of the paper.} Section \ref{s-basic} contains generalities about ARS. Section \ref{s-local} contains the construction (under generic conditions) of the local representations. Section \ref{s-abnormals} contains results on abnormal extremals. Section \ref{s-greg} contains the analysis of the spheres and of the cut locus for the nilpotent points of \Tone. This is the most important and technical part of the paper. Section \ref{s-heat} contains the discussion of the heat diffusion on 3-ARS. Appendix \ref{a-genericity} contains the proof of the 
genericity of (G1), (G2), (G3). Appendix \ref{b-heat} contains the explicit construction of the heat kernel for $\Delta_L$.

\section{Basic Definitions}
\label{s-basic}

 \begin{definition}
 \label{d-ARS}
 A {\it $n$-dimensional almost-Riemannian structure} ($n$-ARS, for short) is a triple
${\mathcal S}=
(E, \f,\langle\cdot,\cdot\rangle)$
 where 
  $E$ is a  vector bundle of rank $n$ over a $n$ dimensional smooth manifold $M$ and $\langle\cdot,\cdot\rangle$ is a Euclidean structure on $E$, that is, $\langle\cdot,\cdot\rangle_q$ is a   scalar product 
on $E_q$ smoothly depending on $q$.  Finally
   $\f:E\rightarrow TM$ is a morphism of vector bundles,    
i.e., {\bf (i)}  the diagram 
\bqn
\label{commutativo}
\xymatrix{
 E  \ar[r]^{\f} \ar[dr]_{\pi_E}   & TM \ar[d]^{\pi}            \\
 & M                          
}    
\eqn
 commutes, where  $\pi:TM\rightarrow M$ and $\pi_E:E\rightarrow M$ denote
  the canonical projections and {\bf (ii)} $\f$ is linear on fibers.

Denoting by $\Gamma(E)$ the $C^\infty(M)$-module of 
smooth sections on $E$, and by 
 $\f_*:\Gamma(E)\rightarrow \VecM$ the map $\sigma\mapsto\f_*(\sigma):=\f\circ\sigma$, we require that  the submodule of Vec$(M)$ given by $\bD=\f_*(\Gamma(E))$ to be bracket generating, i.e.,
  $Lie_q(\bD)= T_qM$ for every $q\in M$.  Moreover, we require that $\f_*$ is injective.

\end{definition}
Here Lie$(\bD)$ is  the smallest Lie subalgebra of Vec(M) containing $\bD$ and 
Lie$_q(\bD)$ is the linear subspace of $T_qM$ whose elements are evaluation at $q$ of elements belonging to Lie$(\bD)$.
The condition that $\bD$ satisfies the Lie bracket generating assumption is also known  as the H\"ormander condition.

We say that a $n$-ARS $(E, \f,\langle\cdot,\cdot\rangle)$  is trivializable if $E$ is isomorphic to the trivial bundle $M 
\times \R^n$. A particular case of  $n$-ARSs is given by $n$-dimensional Riemannian manifolds. In this case $E=TM$ and $\f$ is the identity.

Let ${\mathcal S}=(E,\f,\langle\cdot,\cdot\rangle)$ be a  $n$-ARS on a manifold $M$.
We  denote by $\bD(q)$  the linear subspace $\{V(q)\mid  V\in \bD\}=\f(E_q)\subseteq  T_q M$. 
The set of points in $M$ such that $\dim(\bD(q))<n$ is called {\it singular set} and denoted by $\Zz $.

%
%

If $(\sigma_1,\ldots,\sigma_n)$ is an orthonormal frame for $\langle\cdot,\cdot\rangle$ on
 an open subset $\Omega$ of $M$, an {\it  orthonormal frame} on $\Omega$ is 
 given by  $(\f\circ\sigma_1,\ldots,\f\circ\sigma_n)$.  Orthonormal frames are systems of local generators of $\bD$.

For every $q\in M$  
and every $v\in\bD(q)$ define
$$
\Gq(v)=\min\{\langle u, u\rangle_q \mid u\in E_q,\f(u)=v\}.
$$
For a vector field $X$, we call $\sqrt{\Gq(X(q))}$ the ponctual norm of $X$ at $q$.


An  absolutely continuous curve $\g:[0,T]\to M$  is  admissible for ${\mathcal S}$ 
if   there exists a measurable essentially bounded function 
$$[0,T]\ni t\mapsto u(t)\in E_{\g(t)},
$$ called {\it control function} such that 
$\dot \g(t)=\f(u(t))$  for almost every $t\in[0,T]$. 
Given an admissible 
curve $\g:[0,T]\to M$, the {\it length of $\g$} is  
\bqn
\ell(\g)= \int_{0}^{T} \sqrt{ \gg_{\gamma(t)}(\dot \g(t))}~dt.
\eqnn
The {\it Carnot-Caratheodory distance} (or sub-Riemannian distance) on $M$  associated with 
${\mathcal S}$ is defined as
\bqn\nonumber
d(q_0,q_1)=\inf \{\ell(\g)\mid \g(0)=q_0,\g(T)=q_1, \g\ \mathrm{admissible}\}.
\eqn

{ 
It is a standard fact that $\ell(\g)$ is invariant under reparameterization of the curve $\g$. 
Moreover, if an admissible curve $\g$ minimizes the so-called {\it energy functional} 
$
E(\g)=\int_0^T \gg_{\gamma(t)}(\dot \g(t))~dt
$
with $T$ fixed (and fixed initial and final point)
then $v(t)=\sqrt{\gg_{\gamma(t)}(\dot \g(t))}$ is constant and 
$\g$ is also a minimizer of $\ell(\cdot)$. 
On the other hand, a minimizer $\g$ of $\ell(\cdot)$, such that  $v$ is constant, is a minimizer of $E(\cdot)$ with $T=\ell(\g)/v$.

The finiteness and the continuity of $d(\cdot,\cdot)$ with respect 
to the topology of $M$ are guaranteed by  the Lie bracket generating 
assumption on the $n$-ARS (see \cite{agra-book}).  
The Carnot-Caratheodory distance  endows $M$ with the 
structure of metric space compatible with the topology of $M$ as differential manifold. 

When  the $n$-ARS is trivializable, the problem of finding a curve minimizing the energy between two fixed points  $q_0,q_1\in M$ is naturally formulated as the distributional optimal control problem  with quadratic cost and fixed final time
\bqn
\dot q=\sum_{i=1}^n u_i X_i(q)\,,~~u_i\in\R\,,
~~\int_0^T 
\sum_{i=1}^n u_i^2(t)~dt\to\min,~~q(0)=q_0,~~q(T)=q_1.
\label{eq-localMIN}
\eqn
where $\{X_1,\ldots,X_n\}$ is an orthonormal frame.
}

\subsection{Minimizers and geodesics}
 In this section we briefly recall some facts about 
geodesics in $n$-ARSs. In particular, we define the corresponding exponential map.

 \bdeff A \emph{geodesic} for a $n$-ARS is an admissible curve $\g:[0,T]\to M$ such that
 $ \gg_{\gamma(\cdot)}(\dot \g(\cdot))$
  is constant and, for every sufficiently
small interval $[t_1,t_2]\subset [0,T]$, the restriction
$\g_{|_{[t_1,t_2]}}$ is a minimizer of $\ell(\cdot)$.
A geodesic for which $\gg_{\gamma(\cdot)}(\dot \g(\cdot))=1$ is said to be
parameterized by arclength.

 A $n$-ARS is said to be \emph{complete}
if $(M,d)$ is complete as a metric space.
If the sub-Riemannian
metric is the restriction of a complete Riemannian metric,
then it is complete.

Under the assumption that the manifold is complete, a
version of the Hopf-Rinow theorem (see \cite[Chapter 2]{burago})
implies that the manifold is geodesically complete (i.e. all geodesics are defined for every $t\geq0$) and that for every two points there exists a minimizing geodesic
connecting them.
\edeff
Trajectories minimizing the distance between two points are solutions
of first order necessary conditions for optimality, which in the case
of sub-Riemannian geometry are given by a weak version of the
Pontryagin Maximum Principle (\cite{pontryagin-book}).
\bt\label{t:pmpw}
Let $q(\cdot)$ be a solution of the
minimization problem (\ref{eq-localMIN}) such that
$\gg_{q(\cdot)}(\dot q(\cdot))$  is constant and $u(\cdot)$ be the
corresponding control.
Then there exists a Lipschitz map $p:t\in [0,T] \mapsto
p(t)\in T^{*}_{q(t)}M\setminus\{0\}$  such that one and only one of
the following conditions holds:
\bi
\iii[(i)]
$
\dot{q}=\dfrac{\partial H}{\partial p}, \quad
\dot{p}=-\dfrac{\partial H}{\partial q}, \quad
u_{i}(t)=\la p(t), X_{i}(q(t))\ra,
\\$
where $H(q,p)=\frac{1}{2} \sum_{i=1}^{n} \la p,X_{i}(q)\ra^{2}$.
\vspace{0.2cm}
\iii[(ii)]
$
\dot{q}=\dfrac{\partial {\cal H}}{\partial p}, \quad
\dot{p}=-\dfrac{\partial {\cal H}}{\partial q}, \quad
0=\la p(t), X_{i}(q(t))\ra,
\\$
where ${\cal H}(t,q,p)=\sum_{i=1}^{n}u_{i}(t) \la p,X_{i}(q)\ra$.
\ei
\et
\noindent
For an elementary proof of Theorem \ref{t:pmpw} see \cite{AgrBarBoscbook}.

\brem If $(q(\cdot),p(\cdot))$ is a solution of (i) (resp.\ (ii)) then
it is called a \emph{normal extremal} (resp.\ \emph{abnormal
extremal}). It is well known that if $(q(\cdot),p(\cdot))$ is a normal
extremal then $q(\cdot)$ is a geodesic (see
\cite{AgrBarBoscbook}). This does not hold in general for
abnormal extremals. An admissible trajectory $q(\cdot)$ can be at the
same time normal and abnormal (corresponding to different covectors).
If an admissible trajectory $q(\cdot)$ is normal but not abnormal, we
say that it is \emph{strictly normal}. An abnormal extremal such that $q(\cdot)$ is constant, is called trivial.

\erem


In the following we denote by $(q(t),p(t))=e^{t\vec{H}}(q_{0},p_{0})$
the solution of $(i)$ with initial condition
$(q(0),p(0))=(q_{0},p_{0})$. Moreover we denote by $\pi:T^{*}M\to M$
the canonical projection.

 Normal extremals (starting from $q_{0}$) parametrized by arclength
correspond to initial covectors $p_{0}\in \Lambda_{q_{0}}:=\{p_{0}\in
T^{*}_{q_{0}}M | \, H(q_{0},p_{0})=1/2\}.$
\bdeff  \label{d:exp}
Consider a $n$-ARS. We define the \emph{exponential map} starting from $q_{0}\in M$
as
\bqn \label{eq:expmap}
\EXP_{q_{0}}: \Lambda_{q_{0}}\times \R^{+} \to M, \qquad
\EXP_{q_{0}}(p_{0},t)= \pi(e^{t\vec{H}}(q_{0},p_{0})).
\eqn
\edeff

Notice that each $\EXP_{q_{0}}(p_{0},\cdot)$ is a geodesic. Next, we recall the definition of cut and conjugate time.

\bdeff \label{def:cut} Let $q_{0}\in M$ and
$\g(t)$ an arclength geodesic starting from $q_{0}$.
The \emph{cut time} for $\g$ is $\tcut(\g)=\sup\{t>0,\, \g|_{[0,t]}
\text{ is optimal}\}$. The \emph{cut locus} from $q_{0}$ is the set
$\Cut(q_{0})=\{\g(\tcut(\g)), \g$  arclength geodesic from
$q_{0}\}$.
\edeff


\bdeff \label{def:con} Let $q_{0}\in M$ and
$\g(\cdot)=\EXP_{q_{0}}(p_{0},\cdot)$, $p_0\in\Lambda_{q_{0}}$ a normal arclength geodesic.
The \emph{first conjugate time} of $\g$ is
$\tcon(\g)=\min\{t>0,\  (p_{0},t)$ is a critical point of
$\EXP_{q_{0}}\}$. The \emph{first conjugate locus} from $q_{0}$ is the
set $\Con(q_{0})=\{\g(\tcon(\g)), \g$  normal arclength geodesic from
$q_{0}\}$.

\edeff
It is well known that, for a geodesic $\g$ which is not abnormal, the
cut time $t_{*}=\tcut(\g)$ is either equal to the conjugate time or
there exists another geodesic $\til{\g}$ such that
$\g(t_{*})=\til{\g}(t_{*})$ (see for instance \cite{agrexp}).

Let $ (q(\cdot),p(\cdot))$ be an abnormal extremal. In the following we use the convention that all points of Supp$(q(\cdot))$ are conjugate points.

\brem
In $n$-ARSs,  the exponential map starting from $q_{0}\in \Zz$ is never a local
diffeomorphism in a neighborhood of the point $q_{0}$ itself. As a consequence the
metric balls centered in $q_0$ are never smooth and both the cut and the
conjugate loci from $q_{0}$ are adjacent to the point $q_{0}$ itself
(see \cite{agratorino}).
\erem

To study local properties of $n$-ARSs, it is useful to use local representations.

   \bdeff\label{def:locrep}
A  \emph{local representation} of a $n$-ARS  ${\mathcal S}$ at  a point $q\in M$ is a $n$-tuple of vector fields  $(X_1,\ldots,X_n)$ on $\R^n$ such that there exist: {\bf i) }a neighborhood $U$ of $q$ in $M$, a neighborhood $V$ of the origin in $\R^n$ and a diffeomorphism $\varphi:U\rightarrow V$ such that $\varphi(q)=(0,\ldots,0)$; {\bf ii)} a local orthonormal frame $(F_1,\dots,F_n)$ of ${\mathcal S}$ around $q$, such that  $\varphi_*F_1=X_1$, $\varphi_*F_2=X_2, \ldots, \varphi_*F_n=X_n$ where $\varphi_*$ denotes the push-forward.
\edeff

%
%

\section{Local Representations}
\label{s-local}

The main purpose of this section is to give local representations of 3-ARS under generic conditions.

\bdeff
\label{d-generic}
A property $(P)$ defined for 3-ARSs  is said to be {\it generic}
if for every rank-3 vector bundle $E$   over $M$, $(P)$ holds for every $\f$ in a residual  subset of the set of  morphisms of vector bundles from $E$ to $TM$, 
endowed with the 
$\con^\infty$-Whitney topology, such that $(M,E,\f)$ is a 3-ARS.

\edeff

For the definition of residual subset, see Appendix \ref{a-genericity}.
We have the following
 \bp 
\label{p-generic}
 Consider a 3-ARS. The following conditions are generic for 3-ARSs.

\hspace{.1cm}{(G1)}  dim($\bD(q))\geq 2$
and $\bD(q)+[\bD,\bD](q)=T_qM$  for every $q\in M$;

\hspace{.1cm}{(G2)} $\Zz$ is an embedded two-dimensional 
submanifold of $M$;

\hspace{.1cm}{(G3)}  the points where $\bD(q)=T_q\Zz$ are isolated. 

\ep
For the proof see  Appendix \ref{a-genericity}. In the following we refer to the set of conditions (G1), (G2), (G3) as to the  {\bf (G)} condition.  {\bf (G)}  is the condition under which the main results of the paper are proved.

A way of rephrasing Proposition \ref{p-generic} is the following (see Figure \ref{f-3-punti}): \\[3mm]
{\bf Proposition \ref{p-generic}bis}
{\em Under the condition {\bf (G)} on the 3-ARS there are three types of points:
\bi
\item {\bf Riemannian} points  where $\bD(q)=T_qM$. 
\item {\bf \Tone} points  where $\bD(q)$ has dimension 2 and is transversal to $\Zz$.
\item {\bf \Ttwo} points   where $\bD(q)$ has dimension 2 and is tangent to $\Zz$.
\ei
Moreover \Ttwo\ points are isolated, \Tone\ points form a 2 dimensional manifold and all the other points are Riemannian points.} \\[5mm]
The main result of this section is the following Theorem. 
\bt
\label{t-lokal}
If a 3-ARS satisfies  {\bf (G)}  then
 for every point $q\in M$ there exists a local representation having the form
\bqn
X_1(x,y,z)=\left(
\ba{c} 1\\0\\0\ea
\right),~~~
X_2(x,y,z)=\left(
\ba{c} 0\\\al(x,y,z)\\\beta(x,y,z)\ea
\right),~~~
X_3(x,y,z)=\left(
\ba{c} 0\\0\\\nu(x,y,z)\ea
\right),~~~
\label{eq-local}
\eqn
where $\alpha(0,0,0)=1$, $\beta(0,0,0)=0$. Moreover one of following conditions holds:
\bd
\iii[Riemannian case:] $\nu(0,0,0)=1$.
\iii[\Tone\ case:] 
$\al=1+x \bar \al $, $\beta=x \bar \beta $, $\nu=x \bar \nu$ where $\bar \al$, $\bar \beta$, $\bar \nu$ are smooth functions such that 
$\bar\nu(0,0,z)^2+\bar\beta(0,0,z)^2=1$. Moreover the function $\bar \nu$ may have zeros only on the plane $\{x=0\}$.
\iii[\Ttwo\ case:]  
$\al=1+x \bar \al_1+h\bar \al_2 $, $\beta=x \bar \beta_1+h\bar \beta_2$ and 
$\nu=h\bar\nu$ where $\bar\nu$, $\bar\al_i$ and $\bar \beta_i$, $i=1,2$ are smooth functions  such that $\bar \beta_1(0,0,0)\neq0$, $h(x,y,z)=z-\varphi(x,y)$ with $\varphi$ a smooth function satisfying
$\varphi(0,0)=\frac {\partial\varphi}{\partial x}(0,0)=\frac {\partial\varphi}{\partial y}(0,0)=\frac {\partial^2 \varphi}{\partial x\partial y}(0,0)=0$.
Moreover the function $\bar \nu$ may have zero only on $\{h(x,y,z)=0\}$ which defines a two-dimensional surface.

\ed
\label{t-prenormal}
\et
\brem
Notice that Theorem \ref{t-prenormal} implies that
\bi
\item
 in the \Tone\ case, $\Zz=\{x=0\}$. 

\item in the \Ttwo\ case, the set $\Zz=\{h(x,y,z)=0\}$ has its tangent space at zero equal to span$\{\partial_x,\partial_y\}$. 
Moreover, the set of \Ttwo\ points being discrete, we can assume the generic condition

\hspace{3cm} {\bf (iiibis)} $\bar\nu(0,0,0)\neq0$, $\frac {\partial^2 \varphi}{\partial x^2}(0,0)\neq0$ 
and $\frac {\partial^2 \varphi}{\partial y^2}(0,0)\neq0$.

Condition {\bf (iiibis)} implies that in the \Ttwo\ case $\nu=a z + b x^2+ c y^2+ o(x^2+y^2+|z|)$ where $a$, $b$ and $c$ are not zero.
\ei
\erem

\brem
Notice that in a neighbourhood where the almost-Riemannian structure is expressed in the form (\ref{eq-local}) the corresponding Riemannian metric has the expression (where it is defined)
$$
\left(
\begin{array}{ccc}
 1 & 0 & 0 \\
 0 & \frac{\beta ^2+\nu ^2}{\alpha ^2 \nu ^2} & -\frac{\beta
   }{\alpha  \nu ^2} \\
 0 & -\frac{\beta }{\alpha  \nu ^2} & \frac{1}{\nu ^2} \\
\end{array}
\right).
$$
The corresponding Riemannian volume is (where it is defined)
\bqn
\omega(x)=\frac{1}{\sqrt{\alpha(x)^2\nu(x)^2}}dx\wedge dy\wedge dz
\label{omega}
\eqn
\erem

\subsection{Proof of Theorem \ref{t-prenormal}}

Let $W$ be a 2-dimensional surface transversal to $\distr$ at $q$. 

\bl
There exist a local coordinate system $(x,y,z)$ centered at $q$ such that 
 $W=\{x=0\}$, $t\mapsto(t,y,z)$ is a geodesic transversal to $W$ for every $(y,z)$ 
and the following triple is an orthonormal frame of the metric
$$
X_1(x,y,z)=\left(
\ba{c} 1\\0\\0\ea
\right),~~~
X_2(x,y,z)=\left(
\ba{c} 0\\\al(x,y,z)\\\beta(x,y,z)\ea
\right),~~~
X_3(x,y,z)=\left(
\ba{c} 0\\0\\\nu(x,y,z)\ea
\right),~~~
$$
where $\alpha$, $\beta$ and $\nu$ are smooth functions with $\alpha(0,0,0)\neq 0$.
\el \label{lemma-prenormal}
{\bf Proof} Assume that a coordinate system $(y,z)$ is fixed on $W$ and fix a transversal
orientation along $W$. Then consider the family of geodesics $t\mapsto\gamma_{yz}
(t)$ parameterized by arclength, positively oriented and transversal to $W$ at $(y,z)$.
The map $(x,y,z)\mapsto \gamma_{yz}(x)$ is a local diffeomorphism and hence defines
a coordinate system. In this system $\partial_x(q)$ has norm 1 and is orthogonal, with respect to $\gg_q$, to $T_q\{x=c\}\cap\bD_q$ for any constant $c$ close to 0. Call it $X_1'$. 

Now, since the distribution has dimension at least two at each point, one can find a
 vector field $X_2'$ of the distribution of norm one whose ponctual norm is equal to 1 and 
which is orthogonal to $X_1'$. It is tangent to ${x=c}$ for any constant $c$ close to 0. 
We can fix the coordinate system on $W$ in such a way $X_2'(0,y,z)=\partial_y(0,y,z)$.
If we complete the orthonormal frame with a vector field $X_3'$ we find that 
$$
X_1'(x,y,z)=\left(
\ba{c} 1\\0\\0\ea
\right),~~~
X_2'(x,y,z)=\left(
\ba{c} 0\\ \al_1(x,y,z)\\ \beta_1(x,y,z)\ea
\right),~~~
X_3'(x,y,z)=\left(
\ba{c} 0\\ \mu_1(x,y,z)\\ \nu_1(x,y,z)\ea
\right),
$$
with $\al_1(0,y,z)=1$ and $\beta_1(0,y,z)=0$. Locally $\alpha_1^2(0,0,0)+\mu_1^2(0,0,0)>0$ hence
we can choose the orthonormal frame $(X_1,X_2,X_3)=(X_1', \frac{\alpha_1 X_2'+\mu_1 X_3'}{\sqrt{\alpha_1^2+\mu_1^2}},
\frac{-\mu_1 X_2'+\alpha_1X_3'}{\sqrt{\alpha_1^2+\mu_1^2}})$ satisfying
$$
X_1(x,y,z)=\left(
\ba{c} 1\\0\\0\ea
\right),~~~
X_2(x,y,z)=\left(
\ba{c} 0\\\alpha(x,y,z)\\\beta(x,y,z)\ea
\right),~~~
X_3(x,y,z)=\left(
\ba{c} 0\\0\\\nu(x,y,z)\ea
\right),
$$
where $\alpha =\sqrt{\alpha_1^2+\mu_1^2}$,
 $\beta=\frac{\alpha_1\beta_1+\mu_1\nu_1}{\sqrt{\alpha_1^2+\mu_1^2}}$ and
$\nu=\frac{-\mu_1\beta_1+\alpha_1\nu_1}{\sqrt{\alpha_1^2+\mu_1^2}}$.
\hfill$\blacksquare$

\bigskip

\noindent {\bf End of the proof of Theorem \ref{t-prenormal}:} Let us start with the Riemannian case. In the construction, we are still free in fixing the vertical axis, that is the curve $z\mapsto (0,0,z)$. Let us choose it orthogonal to $X_2$ and parameterized such that $\partial_z(0,0,z)$ has norm one. Then $\mu_1(0,0,z)=0$ and $\nu_1(0,0,z)=1$ for $z$ small enough. As a consequence, since $\alpha_1(0,0,0)=1$ and $\beta_1(0,0,0)=0$ we find at $(0,0,0)$
$$
X_2=\left(
\ba{c} 0\\\sqrt{\alpha_1^2+\mu_1^2}\\\frac{\alpha_1\beta_1+\mu_1\nu_1}{\sqrt{\alpha_1^2+\mu_1^2}}\ea
\right)=\left(\ba{c} 0\\1\\0\ea\right),~~~
X_3=\left(
\ba{c} 0\\0\\\frac{-\mu_1\beta_1+\alpha_1\nu_1}{\sqrt{\alpha_1^2+\mu_1^2}}\ea
\right)=\left(\ba{c} 0\\0\\1\ea\right),
$$
which finishes the proof for the Riemannian case.

Let us assume now that $q\in\Zz$. Since in the construction in the proof of lemma \ref{lemma-prenormal} $X_1$ and $X_2$ are assumed of ponctual norm 1, then 
$X_3$ is zero along $\Zz$. hence $\mu=\nu=0$ on $\Zz$. Hence 
$$
X_2=\left(\ba{c} 0\\ \alpha\\ \beta\ea\right),~~~
X_3=\left(\ba{c} 0\\0\\ \nu\ea\right),
$$
where $\alpha=\alpha_1$, $\beta=\beta_1$ and $\nu=0$ on $\Zz$.

Now, if $\Zz$ is transversal to the distribution, one can fix $W=\Zz=\{x=0\}$
which implies that $\alpha=1$ and $\beta=0$ for $x=0$ since $W=\{x=0\}$, and 
$\nu=0$ for $x=0$ since $\Zz=\{x=0\}$. As a consequence
$\al=1+x \bar \al $, $\beta=x \bar \beta $, $\nu=x \bar \nu$ where $\bar \al$, $\bar \beta$, $\bar \nu$ are smooth functions. The fact that $[\distr,\distr](q)=T_qM$ implies that $\bar\beta(0,0,0)\neq0$ or $\bar \nu(0,0,0)\neq0$. Up to a reparameterization of the $z$-axis, one can hence assume that $\bar\beta(0,0,z)^2+\bar\nu(0,0,z)^2=1$ for $z$ small enough.

Finally, if $\Zz$ is tangent to the distribution then its tangent space at $q$ is generated 
by $X_1(q)$ and $X_2(q)$. This implies that it exists $\varphi$ sur that $\Zz=\{h(x,y,z)=z-\varphi(x,y)=0\}$ where $\varphi$ is a 
smooth function such that 
$\varphi(0,0)=\frac {\partial\varphi}{\partial x}(0,0)=\frac {\partial\varphi}{\partial y}(0,0)=0$. Now, up to a rotation in the 
$(x,y)$-coordinates (and in the choice of $X_1$ and $X_2$) we can moreover assume that 
$\frac {\partial^2 \varphi}{\partial x\partial y}(0,0)=0$. The fact that $\bar\nu$ has no zero outside $\{h(x,y,z)=0\}$ is a
 consequence of the fact that this last set is a two dimensional manifold passing through $q$ which is included in $\Zz$ 
implying that locally $\Zz=\{h(x,y,z)=0\}$.

\section{Abnormal extremals}

\label{s-abnormals}

In this section we investigate the presence and characterization of abnormal extremals for 3-ARS. Notice that, on one hand, there are no abnormal extremals starting from a Riemannian point. On the other hand,  as we will see, from \Tone\ points abnormal extremals  start, except in some exceptional case. Roughly speaking, abnormal extremals can be described as trajectories of a field of directions defined on the surface $\mathcal{Z}$ and corresponding, at a given point $q$ of $\mathcal{Z}$, to the one-dimensional intersection $T_q\mathcal{Z}\cap \blacktriangle(q)$.
Let us formalize this point. 

Assume that {\bf (G)} holds. Let $q_0\in\mathcal{Z}$ be of type-1 and assume that $X_1,X_2,X_3$ are vector fields spanning $\blacktriangle$ in a neighborhood of $q_0$ such that $X_1(q_0)\wedge X_2(q_0)\neq 0$. Assume moreover that $\sum_{i=1}^3 \det\left(X_1(q_0), X_2(q_0),[X_{i-1}, X_{i+1}](q_0)\right)X_i(q_0)\neq 0$, with the convention that $X_0=X_3,\ X_4=X_1$. This condition is satisfied 
 in the whole $\mathcal{Z}$ except some isolated points. Then we claim that there exists an open neighborhood $U$ of $q_0$ such that for any point $q\in U\cap \mathcal{Z}$ there exists only one nontrivial abnormal extremal passing through $q$ and the latter is, up to reparametrization, a trajectory of  the vector field 
\begin{align}
&X(q)=\sum_{i=1}^3 u_i(q)X_i(q)&\label{abn1}\\ 
&u_i(q)=\det\left(X_1(q), X_2(q),[X_{i-1}, X_{i+1}](q)\right),\quad i=1,2,3.\label{abn2}&
\end{align} 



From the Pontryagin maximum principle we have that with  each abnormal extremal $q(\cdot)$ one can associate an adjoint vector $p(\cdot)$ such that $\langle p(t),X_i(q(t))\rangle=0$, for $i=1,2,3$. It turns out that $q(\cdot)$ must be contained in $\mathcal{Z}$ and, in a neighborhood of $q_0$, $p(t)$ must be proportional to the nonzero vector $X_1(q(t))\wedge X_2(q(t))$. 
By differentiating with respect to time the equality $\langle p(t),X_i(q(t))\rangle=0$, and knowing that the adjoint vector $p$ satisfies the equation 
\[\dot p=\sum_{i=1}^3 u_i\left(\frac{\partial X_i}{dq}\right)^Tp\] 
we get that
$\sum_{j\neq i} u_j \langle p,[X_i,X_j]\rangle=0$
leading to 
\[\sum_{j\neq i} u_j \langle X_1\wedge X_2,[X_i,X_j]\rangle=\sum_{j\neq i} u_j\, \det( X_1, X_2,[X_i,X_j])=0.\]
Whenever the triple of components $\mathrm{det}( X_1, X_2,[X_i,X_j])$, $i\neq j$ is different from $0$ we get that the linear equation above is satisfied for the $u_i$'s given in \eqref{abn2}.  Taking into account the local representation given in Theorem~\ref{t-prenormal} one can see by a direct computation, and by using the fact that $\partial_{x} \beta(0,0,0)\neq 0$ or  $\partial_{x} \nu(0,0,0)\neq 0$, that such a triple is always nonzero for a type-1 point. Note that the condition $X(0)\neq 0$ characterizing the possibility of having a nontrivial abnormal extremal parameterized by arclength passing through the origin  is verified whenever $\bar\nu(0,0,0)\neq 0$.
Note also that, close to a type-1 point, the equations $\langle p,X_i(q)\rangle=0,\ i=1,2,3,$ define a three-dimensional submanifold of $T^*M$ (this can be checked easily via the local representation of $X_1,X_2,X_3$). The Hamiltonian field, with $u_i=u_i(q)$, turns out to be tangent to such submanifold, confirming that the abnormal extremals are exactly those trajectories that satisfy \eqref{abn1}-\eqref{abn2}. 


On a type-2 point, again by direct computation with the local representation defined as in Theorem~\ref{t-prenormal}, one sees that the vector field $X$ vanishes at $q_0=0$. By considering $x,y$ as local coordinates in $\mathcal{Z}$, we have the following linearized equation for the abnormal extremals around the type-2 point
\begin{eqnarray*}
\dot x & = & 2b y-x\bar{\beta}_1(0,0,0)\\
\dot y & = & -2a x\\
\end{eqnarray*}
where $a=\frac{\partial^2 \varphi}{\partial x^2}$ and $b=\frac{\partial^2 \varphi}{\partial y^2}$. Note that, depending on the values $a,b,\bar{\beta}_1(0,0,0)$, the previous system can be stable or  unstable, and may have real or complex non-real eigenvalues. Moreover, since $u_1(0)=-\bar{\beta}_1(0,0,0)\neq 0$, it turns out that abnormal extremals parameterized by arclength cannot reach or escape from a type-2 point in finite time.

Concerning optimality of abnormal extremals, since $\bD(q)+[\bD,\bD](q)=T_qM$ for every $q\in M$, as a consequence of  a Theorem of A.Agrachev and J.P. Gauthier  \cite{AgrBarBoscbook,agra-gau-premio} if an abnormal extremal is optimal then it is not strict. Generically this can never happen. 
Let us notice that, since the vector field $X_3$ is zero on $\mathcal{Z}$ with respect to the chosen local representation, optimality would imply that $u_3=0$ locally along the trajectory, and this implies $\partial_x\beta=\bar{\beta}=0$ along the trajectory.

%
%
%

\begin{figure}
\begin{center}
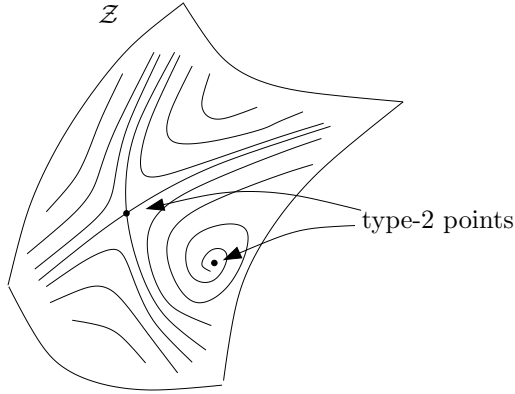
\caption{The singular set $\Zz$ with the field of abnormal extremals. The \Ttwo\ points correspond to to singularity of this field. All other points on $\Zz$ are \Tone\ points.\label{abnormal-extremals-3D}}
\end{center}
\end{figure}

The following theorem summarizes the results obtained in this section.
\begin{theorem}
For any type-1 point there exists an abnormal extremal, parameterized by arclength passing through it
if and only if
\[\sum_{i=1}^3 \det\left(X_1(q_0), X_2(q_0),[X_{i-1}, X_{i+1}](q_0)\right)X_i(q_0)\neq 0,\] with the convention that $X_0=X_3,\ X_4=X_1$. Using the local representation given by Theorem~\ref{t-prenormal}, this condition is equivalent to $\bar\nu(0,0,0)\neq 0$. There is no nontrivial abnormal extremal passing through type-2 points, which are poles of the extremal flow corresponding to abnormal extremals.

Generically all nontrivial abnormal extremals  are not optimal and they are strictly abnormal.
\end{theorem}

\brem
If $q_0$ is a  \Tone\ point such that $\sum_{i=1}^3 \det\left(X_1(q_0), X_2(q_0),[X_{i-1}, X_{i+1}](q_0)\right)X_i(q_0)=0,$ then the trajectory $q(\cdot)=q_0$ is a trivial abnormal extremal.
\erem

\section{Nilpotent approximations}
\label{s-greg}
For each kind of points it is an easy exercise to find the nilpotent approximation in the coordinate system constructed in the local representation. 
For the general theory of the nilpotent approximation, see, for instance, \cite{AgrBarBoscbook,bellaiche}. We have:

\bi
\item in the Riemannian case 
$$
\hat X_1(x,y,z)=\left(
\ba{c} 1\\0\\0\ea
\right),~~~
\hat X_2(x,y,z)=\left(
\ba{c} 0\\1\\0\ea
\right),~~~
\hat X_3(x,y,z)=\left(
\ba{c} 0\\0\\1\ea
\right),
$$

\item in the \Tone\ case
$$
\hat X_1(x,y,z)=\left(
\ba{c} 1\\0\\0\ea
\right),~~~
\hat X_2(x,y,z)=\left(
\ba{c} 0\\1\\\cos(\sigma)x\ea
\right),~~~
\hat X_3(x,y,z)=\left(
\ba{c} 0\\0\\\sin(\sigma)x\ea
\right),
$$
where $\sigma\in[0,\pi/2]$ is a parameter.

\item in the \Ttwo\ case  

$$
\hat X_1(x,y,z)=\left(
\ba{c} 1\\0\\0\ea
\right),~~~
\hat X_2(x,y,z)=\left(
\ba{c} 0\\1\\x\ea
\right),~~~
\hat X_3(x,y,z)=\left(
\ba{c} 0\\0\\0\ea
\right).
$$
\ei
\brem
For \Tone\ points, the nilpotent approximation is not universal and keeps track of the original vector fields through the parameter $\sigma$. The nilpotent approximation for the \Ttwo case is a special case of the one obtained for the \Tone\ case.
\erem

The computation of the exponential flow of the nilpotent approximation is trivial in the Riemannian case.
In the \Tone\ case (including also the \Ttwo\ case), the Hamiltonian 
for the normal flow is given by
$$
H(p,q)=\frac12(p_x^2+(p_y+x\cos(\sigma)p_z)^2+x^2\sin(\sigma)^2p_z^2).
$$
One computes easily that the geodesic with initial conditions $x(0)=y(0)=z(0)=0$, $p_x(0)=\cos(\theta)$, 
$p_y(0)=\sin(\theta)$, $p_z(0)=a$ is given, when $a\neq 0$, by
\begin{eqnarray*}
x(a,\theta,t) &= & \frac 1 a (\cos (\theta)\sin (a t) + (\cos (a t) - 
       1)\cos (\sigma)\sin (\theta)),\\
y(a,\theta,t) &= & \frac 1 a ((\sin (a t) - a t) \sin (\theta ) \cos ^2 (\sigma ) -
 (\cos (a t) - 1) \cos (\theta ) \cos (\sigma ) + a t \sin (\theta )),\\
z(a,\theta,t) &= & \frac 1{8 a^2} \bigg(
4\sin(2\theta)\cos(\sigma)\cos(at)(1-\cos(at))\\
&&+\cos(2\theta)\big(2at\sin^2(\sigma)-\sin(2at)(1+\cos^2(\sigma))+4\sin(at)\cos^2(\sigma)\big)\\
&&+2at(1+\cos^2(\sigma))-\sin(2at)\sin^2(\sigma)-4\sin(at)\cos^2(\sigma)\bigg)
\end{eqnarray*}
and, when $a=0$, by
\begin{eqnarray*}
x(a,\theta,t) &= & t \cos (\theta),\\
y(a,\theta,t) &= & t \sin (\theta ),\\
z(a,\theta,t) &= & \frac 1 4 t^2\cos(\sigma)\sin(2\theta).
\end{eqnarray*}

\noindent {\bf Notations.}
We denote $\gamma(a,\theta,t)=(x(a,\theta,t),y(a,\theta,t),z(a,\theta,t)).$
In the following, we denote $\tau$ the smallest positive real number such that 
\begin{equation}
\sin (\tau) \cos ^2(\sigma )+ \tau \cos (\tau) \sin ^2(\sigma )=0.\label{cut-time}
\end{equation}
The number $\tau$ belongs to $[\frac\pi 2, \pi]$. We also denote by $s_1$ the first positive real number such that 
\begin{equation}
s_1=\tan(s_1), \label{conjugate-time}
\end{equation}
and by $\thep$ and $\them$ the angles defined modulo $\pi$ such that
\begin{eqnarray}
\cos (\sigma ) \cos (\thep ) \sin(\tau)+\cos(\tau) \sin (\thep)&=&0,\label{angle-cc+}\\
-\cos (\sigma ) \cos (\them) \sin(\tau)+\cos(\tau) \sin (\them)&=&0.\label{angle-cc-}
\end{eqnarray}
One checks easily that $\them=-\thep$.

\subsection{Conjugate time in the nilpotent cases}

In this section we prove the following result.

\bt
Let us consider a type-1 nilpotent point for a fixed  value of $\sigma$
then 
\bi
\item if $\cos(\sigma)=0$ any geodesic with initial conditions $x(0)=y(0)=z(0)=0$, $p_x(0)=\cos(\theta)\neq 0$,  $p_y(0)=\sin(\theta)$ and $p_z(0)=a\neq 0$
has a first conjugate time equal to $\frac{s_1}{|a|}$. If $a=0$ and $\cos(\theta)\neq 0$, the geodesic has no conjugate time. If $\cos(\theta)=0$ then the geodesic is entirely included in the conjugate locus.
\item if $\cos(\sigma)\neq 0$ 
any geodesic with initial conditions $x(0)=y(0)=z(0)=0$, $p_x(0)=\cos(\theta)$,  $p_y(0)=\sin(\theta)$ and $p_z(0)=a\neq 0$
has a first conjugate time in the interval $[\frac{2\tau}{|a|},\frac{2\pi}{|a|}]$.
If $a=0$, the geodesic has no conjugate time.
\ei
\et
{\bf Proof.} In the following, instead of considering the case $a<0$ we equivalently consider the case $a>0$ with $t<0$.

The computation of the Jacobian of the exponential map gives, for $a\neq 0$
$$
Jac=\det(\frac{\partial \gamma}{\partial a},\frac{\partial \gamma}{\partial \theta},\frac{\partial \gamma}{\partial t}) =
\frac{1}{2a^4}(A\cos(2\theta)+B+C\sin(2\theta))
$$
with
\begin{eqnarray*}
A&=&at\sin ^2(\sigma ) (at \cos (at)-\sin (at)), \\
C&=&-at \cos (\sigma )\sin ^2(\sigma ) (2 \cos (at)+at\sin (at)-2), \\
B&=&4 (\cos(at)-1) \cos ^2(\sigma )+\frac{at}{2} \left(2 at \cos (at) \sin ^2(\sigma )+3\cos(2 \sigma) \sin(at)+\sin(at)\right),
\end{eqnarray*}
and, for $a=0$
$$
Jac=-\frac{1}{12} t^4 (1+\sin^2(\sigma)(1+2\cos(2\theta))).
$$

For $a=0$ and $t\neq 0$, one can check easily that $Jac=0$ if and only if $\sigma$ and $\theta$ are equal to $\frac\pi 2 [\pi]$. This allows to prove the cases corresponding to $a=0$.

Assume $a> 0$. If $t$ is fixed, there exists a conjugate point for a certain $\theta$ if and only if $B^2-A^2-C^2\leq0$.
After simplification, one gets
\begin{eqnarray*}
B^2-A^2-C^2&=&64 \cos ^2(\sigma )\sin \left(\frac{at}{2}\right) \left(\frac{at}{2} \cos \left(\frac{at}{2}\right)- \sin \left(\frac{at}{2}\right)\right)  \\
&&\times\left(
   \sin \left(\frac{at}{2}\right) \cos ^2(\sigma )+ \frac{at}{2} \cos \left(\frac{at}{2}\right) \sin ^2(\sigma )\right)\\
&&\times \left( \frac{at}{2} \cos
   \left(\frac{at}{2}\right)- \sin \left(\frac{at}{2}\right)- \left(\frac{at}{2}\right)^2 \sin \left(\frac{at}{2}\right) \sin ^2(\sigma )\right).
\end{eqnarray*}

The term $\sin \left(\frac{at}{2}\right)$ is positive if $0<at<2\pi$.
The term $\left(\frac{at}{2} \cos \left(\frac{at}{2}\right)- \sin \left(\frac{at}{2}\right)\right)$ 
is negative if $0<at<2\pi$.
The last term is negative for $0<at<2\pi$, being the sum of 
$\left(\frac{at}{2} \cos \left(\frac{at}{2}\right)- \sin \left(\frac{at}{2}\right)\right)$
and 
$- \left(\frac{at}{2}\right)^2 \sin \left(\frac{at}{2}\right) \sin ^2(\sigma )$ 
which are both negative for $0<at<2\pi$. 
Since $\tau$, defined by (\ref{cut-time}) belongs to $[\frac\pi 2, \pi]$, the smallest time $t_1(a)$ such that 
$\sin \left(\frac{at}{2}\right) \cos ^2(\sigma )+ \frac{at}{2} \cos \left(\frac{at}{2}\right) \sin ^2(\sigma )=0$ 
belongs to $[\frac\pi a, \frac{2\pi} a]$. 

The same computations can be done with $t<0$. In that case $t_1(a)=-\frac{2\tau}{a}$ and belongs to $[-\frac{2\pi} a, -\frac\pi a]$. 

If $\cos(\sigma)\neq 0$ and $t>0$, $t_1(a)$ is the first time for which $B^2-A^2-C^2\leq 0$ and, as a consequence, for any $\theta$ the geodesic with initial data $(a,\theta)$ is not conjugate at time $t<t_1(a)$. At time $t_1(a)$, since $B^2-A^2-C^2=0$ there are exactly
two values of $\theta$ in $[0,2\pi[$ such that the jacobian is zero, when just after time $t_1(a)$ there are 4. 
One can check easily that
if $t=\frac{2\pi}{a}$ then $A=B=4\pi^2\sin^2(\sigma)$ and $C=0$ which implies that $Jac\geq 0$ for any $\theta$. Moreover for $0<t<t_1(a)$ we know that the jacobian is not zero. But 
$B=\frac{1}{12} (at)^4 (\cos (2 \sigma)-3)+o((at)^5)$ which is negative hence for $t$ small $Jac< 0$.
Hence we know that the conjugate time of the geodesic is between $t_1(a)$ and $\frac{2\pi}{a}$.
The same arguments work for the part of the synthesis corresponding to $t<0$.

If $\cos(\sigma)= 0$, then $B^2-A^2-C^2= 0$ for all $t$. It corresponds to the fact that, in that case, for every $a>0$ and $t>0$, 
$\frac{\partial \gamma}{\partial a}(a,\frac\pi 2,t)=\frac{\partial \gamma}{\partial a}(a,-\frac\pi 2,t)=0$. The jacobian is equal to
$$
Jac=\frac{t\cos(\theta)^2}{a^3}(at\cos(at)-\sin(at)).
$$
When $\theta\neq \frac \pi 2 [\pi]$ the first conjugate time is $t=\frac{s_1}{a}$ where $s_1$ is defined by (\ref{conjugate-time}).

\hfill$\blacksquare$

\begin{figure}
\begin{center}
\includegraphics[width=5cm]{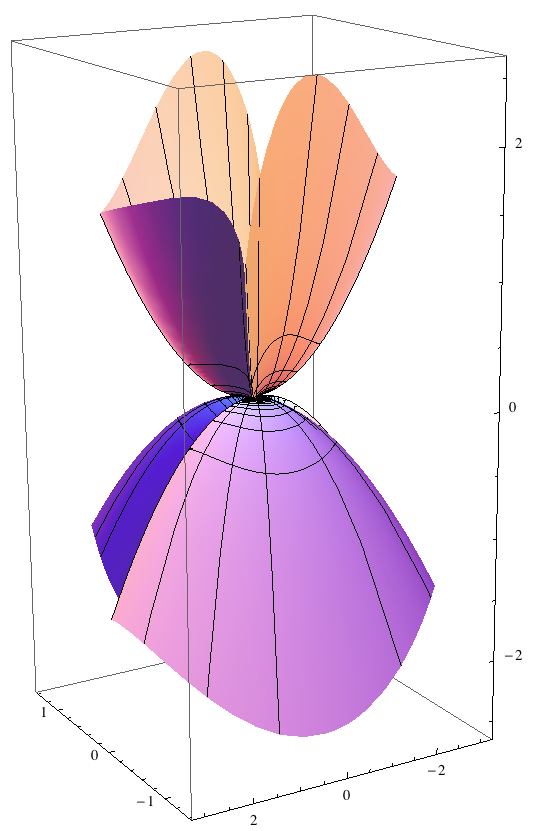}
\caption{The first conjugate locus in the case $\sigma=\pi/4$\label{f-conj-1}}
\end{center}
\end{figure}

\subsubsection{Some numerical simulations describing the conjugate locus}
One can go further in the study  of the conjugate locus in the nilpotent case. This is out of the purpose of this paper. Let us just mention that the first conjugate locus for $\sigma\in]0,\pi/2[$ looks like a suspension of a 4-cusp astroid, similarly to the 3D contact case. The interesting difference is that the two components of the conjugate locus for $z>0$ and $z<0$ are  twisted of an angle which depends on $\sigma$, see Figure \ref{f-conj-1}.

\subsection{Cut locus in the nilpotent cases}

In this section we prove the following result.
\bt
Let us consider a type-1 nilpotent case for a fixed value of $\sigma$.
\bi
\item Case $\cos(\sigma)=0$.
For an initial condition with $\cos(\theta)= 0$ or $a=0$, the corresponding geodesic is optimal for any time.  
For an initial condition with $\cos(\theta)\neq 0$ and $a\neq 0$, the cut time is equal to $\frac{\pi}{|a|}$.
The cut locus at $(0,0,0)$ is the set 
$$\{(x,y,z)\; | \; x=0, z\neq 0\}.$$
\item Case $\cos(\sigma)\neq 0$. 
For an initial condition with $a=0$, the corresponding geodesic is optimal for any time. 
For an initial condition with $a\neq 0$, the cut time is equal to $\frac{2\tau}{|a|}$ with $\tau$ 
defined by (\ref{cut-time}).
The cut locus at $(0,0,0)$ is 
$$
\bigcup_{a\neq0,\theta\in\R}\gamma(a,\theta,\frac{2\tau}{|a|}).
$$ 
More precisely the cut locus is included in the union of the half plane $P_+=\{z>0,\; \cos(\thep)y+\sin(\thep)x=0\}$ 
and the half plane $P_-=\{z<0,\; \cos(\them)y+\sin(\them)x=0\}$ where $\thep$ and $\them$ are defined by 
(\ref{angle-cc+}) and (\ref{angle-cc-}). The intersection of the cut locus with $P_+$ is exactly the set of points which are 
on or above the curve 
$$
a\mapsto (x_\sigma^+(a),y_\sigma^+(a),z_\sigma^+(a))
$$
and the intersection of the cut locus with $P_-$ is exactly the set of points which are 
on or below the curve 
$$
a\mapsto (x_\sigma^-(a),y_\sigma^-(a),z_\sigma^-(a))
$$
where
\begin{eqnarray}
x_\sigma^\pm(a) &= &\frac 2 a \tan(\pm\tau)(\cos^2(\tau)+\sin^2(\tau)\cos^2(\sigma))\cos(\theta_\sigma^\pm),\\
y_\sigma^\pm(a)&=& -\frac 2 a \tan(\pm\tau)(\cos^2(\tau)+\sin^2(\tau)\cos^2(\sigma))\sin(\theta_\sigma^\pm),\\
z_\sigma^\pm(a)&=&\frac{\pm\tau}{a^2}-\frac{\tan(\pm\tau)}{a^2}(\cos^2(\tau)-\sin^2(\tau))(\cos^2(\tau)+\sin^2(\tau)\cos^2(\sigma)).
\end{eqnarray}
\ei
\et

\subsubsection{Case $\cos(\sigma)\neq 0$}

We make the proof for the upper part of the cut locus, the computations being the same for the lower part.

Let us make the following observation: if one consider the closed curve $\theta\mapsto(x,y)(a,\theta,t)$, it happens to be an ellipse for any value of $a$ and $t$. The ellipse is flat when the coefficients of $\cos(\theta)$ and $\sin(\theta)$ in $x$ and $y$ form a matrix of zero determinant which gives the equation
\begin{equation}
\sin \left(\frac{at}{2}\right)\left(\sin \left(\frac{at}{2}\right) \cos ^2(\sigma)+\frac{at}{2}\cos \left(\frac{at}{2}\right) \sin ^2(\sigma)\right)=0.
\end{equation} 
One easily proves that the first positive time satisfying this relation is $t_1(a)=\frac{2\tau}{a}$ computed before, where $\tau$ is defined
by (\ref{cut-time}).
Along this flat ellipse, the extremities correspond to values of $\theta$ such that $\partial_\theta x=\partial_\theta y=0$.
Since at $t=t_1(a)$
$$\partial_\theta x=
-\frac{2}a \sin (\tau) \left(\cos (\sigma ) \cos (\theta ) \sin (\tau)+\cos(\tau) \sin (\theta )\right)
$$
then the values of $\theta$ corresponding to the extremities are the solutions $\thep$ of (\ref{angle-cc+}).
In particular it does not depend on $a$.
Thanks to the fact that for $t=t_1(a)$ the curve $\theta\mapsto(x,y)(a,\theta,t_1(a))$ is a flat ellipse, one gets
$x(a,\thep+\vartheta,t_1(a))=x(a,\thep-\vartheta,t_1(a))$ and $y(a,\thep+\vartheta,t_1(a))=y(a,\thep-\vartheta,t_1(a))$.

Let us consider now the variable $z$. Its derivate with respect to $\theta$ satisfies
\begin{eqnarray*}
\partial_{\theta}z&=&\frac{1}{4a^2}(8 \cos (at) \cos (\sigma)\sin ^2\left(\frac{at}{2}\right) \cos (2 \theta ) \\
&&+\left((2 at-4 \sin (at)+\sin (2 at)) \cos ^2(\sigma)-2at+\sin (2 at)\right) \sin (2 \theta ))
\end{eqnarray*}
hence it is a function of $\theta$ of the type $A_{t,a}\cos(2\theta+B_{t,a})$. As a consequence, if $\theta_0$ is a zero of $\partial_{\theta}z$
then $z(\theta_0-\vartheta)=z(\theta_0+\vartheta)$. Now fixe $t=t_1(a)$ and $\theta=\thep$.
Combining equation (\ref{cut-time}) and (\ref{angle-cc+}), one prove that $(\cos(\thep),\sin(\thep))$ is colinear to 
$(\cos(\tau), -\cos(\sigma) \sin(\tau))$ which implies that $(\cos(2\thep),\sin(2\thep))$ is colinear to 
$(\cos(\tau)^2 - \cos(\sigma)^2 \sin(\tau)^2 , -2 \cos(\sigma) cos(\tau) \sin(\tau))$. Replacing 
in the formula of $\partial_{\theta}z$ one finds
$$
\partial_{\theta}z(a,\thep,t_1(a))=\frac{\lambda}{a^2} (\cos(\sigma) \sin(\tau) (2 \cos(\sigma)^2 \sin(\tau) + 2\tau \cos(\tau) \sin(\sigma)^2)=0
$$
thanks to equation (\ref{cut-time}).
This proves that $\thep$ is such that $z(a,\thep+\vartheta,t_1(a))=z(a,\thep-\vartheta,t_1(a))$. But we have yet proved that
$(x,y)(a,\thep+\vartheta,t_1)=(x,y)(a,\thep-\vartheta,t_1)$.
Hence we have proved that the lift of the flat ellipse (except its extremities) is included in the Maxwell set of points where two geodesics
of same length intersect one each other. 
Moreover for what concerns the two points correponding to the extremities of the flat ellipse, since $\partial _\theta\gamma (a,\theta_1,t_1(a))=0$, they are in the conjugate locus.


In what follows, we prove that the union for $a> 0$ of the flat ellipses corresponding to $t=t_1(a)$ is in fact the upper part of the cut locus. 

Let us first prove that the ellipses corresponding to $(a,t)$ with $0\leq a<2\tau$ 
and $t=1$ have no intersection. In order to prove that they are disjoint, we are going to prove 
that their projections on the $(x,y)$-plane are disjoint. We compute
the determinant 
$$
{\mathfrak D}=\left|
\ba{cc}
\partial_a x & \partial_\theta x\\
\partial_a y & \partial_\theta y
\ea
\right|.
$$
If we prove that it is never zero for every $a$ smallest that the one corresponding to the flat ellipse that is $2\tau$, then it is of constant sign proving that the vector $(\partial_a x(a,\theta,1), \partial_a y(a,\theta,1))$
points inside the ellipse for every $0\leq a<2\tau$ and every $\theta$. As a consequence we get that, before the flat ellipse which is singular, all the ellipses
are disjoint. The computation gives :
$$
{\mathfrak D}={\mathfrak A}\cos(2\theta)+{\mathfrak B}+{\mathfrak C} \sin(2\theta)
$$
where
\begin{eqnarray*}
{\mathfrak A} &= & \frac{1}{2} a (a \cos (a)-\sin (a)) \sin ^2(\sigma),\\
{\mathfrak B} &= & \frac{1}{2} \left(a (a \cos (a)-\sin (a))-\cos ^2(\sigma) \left(\left(a^2-4\right) \cos (a)-3 a \sin (a)+4\right)\right),\\
{\mathfrak C} &= & -a \cos (\sigma) \left(a \cos \left(\frac{a}{2}\right)-2 \sin \left(\frac{a}{2}\right)\right) \sin \left(\frac{a}{2}\right) \sin
   ^2(\sigma).
\end{eqnarray*}
If ${\mathfrak A} ^2+{\mathfrak C}^2-{\mathfrak B}^2<0$, then ${\mathfrak D}$ has the same signe 
as ${\mathfrak B}$ whatever $\theta$. But
\begin{eqnarray*}
{\mathfrak A} ^2+{\mathfrak C}^2-{\mathfrak B}^2&=&-\frac{1}{2} \cos ^2(\sigma) \left(a \cos \left(\frac{a}{2}\right)-2 \sin \left(\frac{a}{2}\right)\right) \sin
   \left(\frac{a}{2}\right) \\
&&\times\left(\sin \left(\frac{a}{2}\right) \cos ^2(\sigma)+\frac{1}{2} a \cos \left(\frac{a}{2}\right) \sin
   ^2(\sigma)\right)\\
&&\times \left(16 \left(\frac{1}{2} a \cos \left(\frac{a}{2}\right)-\sin \left(\frac{a}{2}\right)\right)-4 a^2 \sin
   \left(\frac{a}{2}\right) \sin ^2(\sigma)\right).
\end{eqnarray*}
The term $\sin\left(\frac{a}{2}\right)$ is positive for $a\in]0,2\tau[$ since $2\tau\leq2\pi$.
The term $\left(a \cos \left(\frac{a}{2}\right)-2 \sin \left(\frac{a}{2}\right)\right)$ is negative for $a\in]0,2\pi[$,
the positive solution of $s \cos(s)-\sin(s)=0$ being greater then $\pi$. The last factor is also negative for $a\in]0,2\pi[$
being the sum of two negative terms on this interval. The remaining factor is positive for $a\in]0,2\tau[$ since it is the one
defining $\tau$ in (\ref{cut-time}). As a consequence ${\mathfrak D}$ is negative for $a\in]0,2\tau[$ whatever $\theta$
and we can conclude that no couple of geodesics of length 1 with $0\leq a,a'<2\tau$ do intersect at time 1. 


A geodesic with the initial condition $(a',\theta)$ with $a'>2\tau$ and $\theta\neq \thep$ is not 
optimal at time 1 since it joins the Maxwell set at time $\frac{2\tau}{a'}<1$.

For what concerns a geodesic with initial condition $(a',\thep)$,
it is not optimal after time $\frac{2\tau}{a'}$. This is due to the fact that it is a strictly normal geodesic which implies that it 
is not optimal after the first conjugate time.

Now consider two geodesics corresponding to $(a',\theta')$ and $(a",\theta")$ with $a'$ and $a"$ less or equal to $2\tau$. 
Let $t_2<1\leq\max(t_1(a'),t_1(a"))$. Reproducing the argument we have developped for $t=1$ we can deduce that
$(x,y)(a',\theta',t_2)\neq (x,y)(a",\theta",t_2)$ which implies that these two geodesics do not intersect at any time $t_2<1$. 

To conclude, we have proved that the sphere of radius 1 
is given by the union of the lifts of the ellipses for $-2\tau\leq a\leq 2\tau$, and that the upper part of the cut locus is exactly the union
for $a> 0$ of the lifts of the flat ellipses corresponding to $t=t_1(a)$.

For what concerns the expressions given in the theorem for $x_\sigma^+$, $y_\sigma^+$, etc, it is just a matter of making 
simplifications in the expression of $\gamma(a,\thep,t_1(a))$ using (\ref{cut-time}) and (\ref{angle-cc+}).

\subsubsection{Case $\cos(\sigma)=0$}

In that case 
\begin{eqnarray*}
x(a,\theta,t) &= & \frac 1 a \cos(\theta)\sin (a t),\\
y(a,\theta,t) &= & t \sin(\theta),\\
z(a,\theta,t) &= & \frac 1{4 a^2} \cos(\theta)^2(2at-\sin(2at)),
\end{eqnarray*}
and 
\begin{eqnarray*}
Jac&=&\frac{t\cos(\theta)^2}{a^3}(at\cos(at)-\sin(at)).
\end{eqnarray*}
Let us again fix $t=1$. For a given $a$, the curve $\theta\mapsto(x(a,\theta,1),y(a,\theta,1))$ is an ellipse.
For $a=\pi$ the ellipse is flat and $(x,y,z)(\pi,\theta,1)=(x,y,z)(\pi,\pi-\theta,1)$. This implies that a geodesic
with initial condition $(a,\theta)$ is no more optimal after time $t=\frac{\pi}{a}$ if $\theta\neq \frac{\pi}{2}[\pi]$.

For what concerns the geodesics with initial condition $\theta=\frac\pi2[\pi]$, one proves easily that they are optimal
for every $t$. It is a simple consequence of the fact that the projection of a curve on the $(x,y)$-plane with the Euclidean
metric preserves its length and that the geodesics with $\theta=\frac\pi2[\pi]$ are geodesics for this last metric.
Moreover, as seen before, they are entirely conjugate.

The ellipses $\theta\mapsto(x(a,\theta,1),y(a,\theta,1))$
with $0\leq a<\pi$ have exactly two common points : $(0,-1,0)$ and $(0,1,0)$. If we consider these ellipses without these two
points, they are disjoint. The same arguments as before allow to conclude that the sphere of radius $t>0$ is the union of the 
lifts of the ellipses with $-\frac\pi t\leq a\leq\frac\pi t$ and that the cut locus is the set $\{(x,y,z)\; | \; x=0, z\neq 0\}$.

\noindent{\bf Remark.} A consequence of the previous computations is that the spheres of the nilpotent cases are sub-analytic.

\subsection{Images of the balls in the nilpotent cases}
In the Riemannian case, the balls are the one of the Euclidean case. 
In the \Ttwo\ case, $\hat X_3$ being null and the couple 
$(\hat X_1, \hat X_2)$ being one representation of the Heisenberg metric,
the balls are those of the Heisenberg case in the corresponding representation.
For what concerns the \Tone\ case, the nilpotent approximation has 
a parameter $\sigma$ and the balls vary with the $\sigma$.
\begin{figure}[here!]
\begin{center}
\includegraphics[width=4truecm]{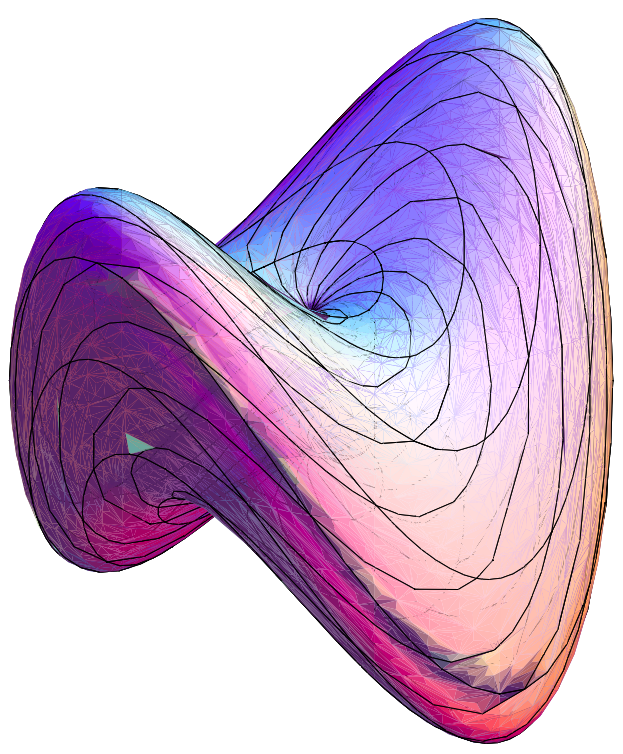}
~~~~~\includegraphics[width=6truecm]{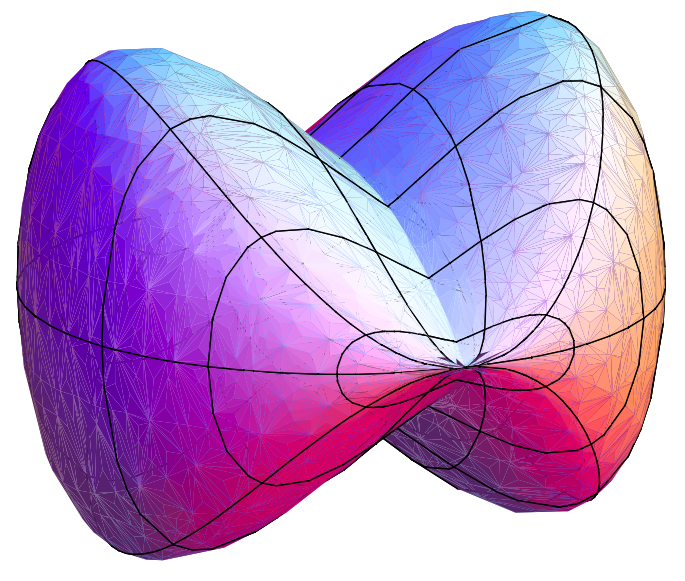}
\caption{The spheres in the case $\sigma=0$ (Heisenberg) and $\sigma=\pi/2$ (Baouendi-Goulaouic)}\label{fig-heisenberg}
\end{center}
\end{figure}
\begin{figure}[here!]
\begin{center}
\includegraphics[width=6truecm]{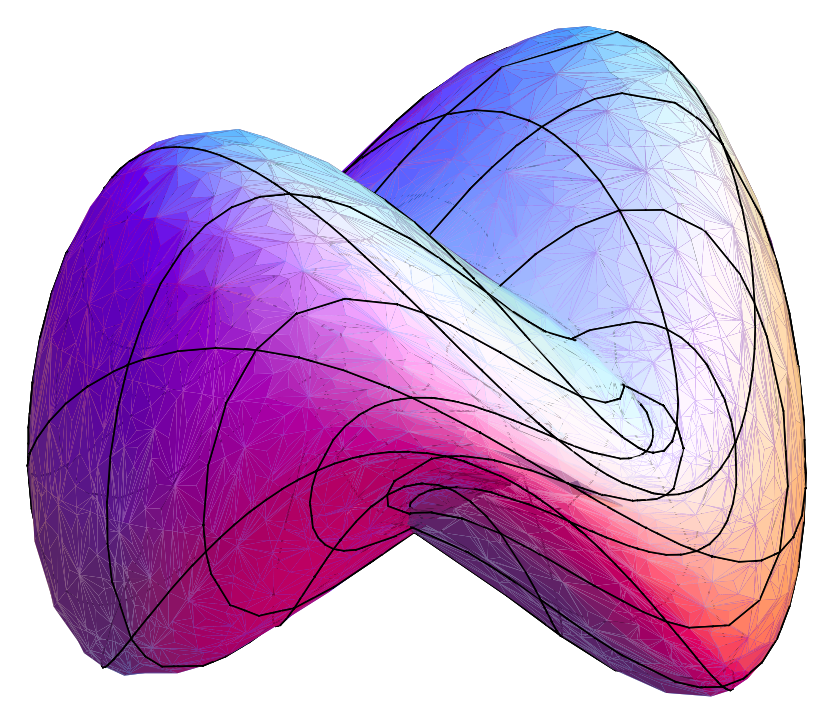}
~~~~~\includegraphics[width=6truecm]{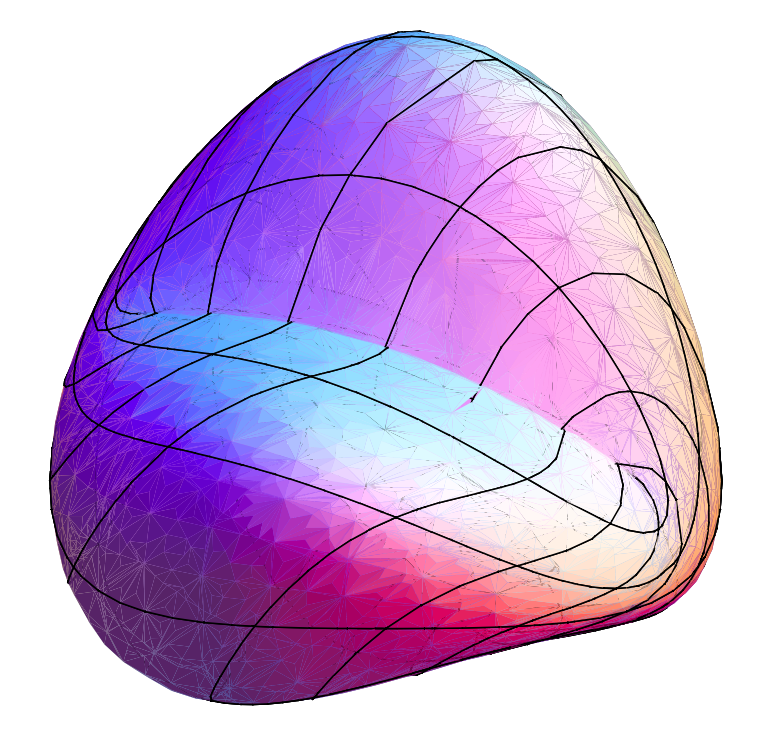}
\caption{Two points of view of the case $\sigma=1$}\label{fig-mixed}
\end{center}
\end{figure}
\begin{figure}[here!]
\begin{center}
\includegraphics[width=6truecm]{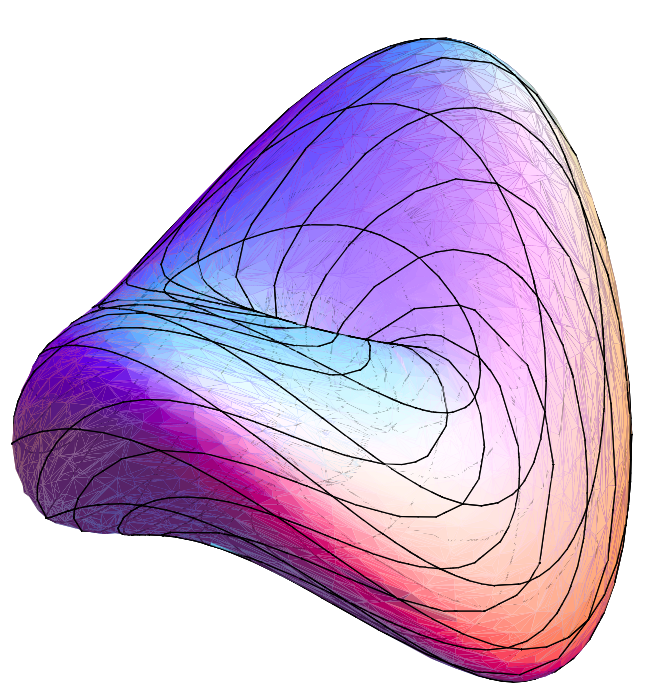}
~~~~~\includegraphics[width=6truecm]{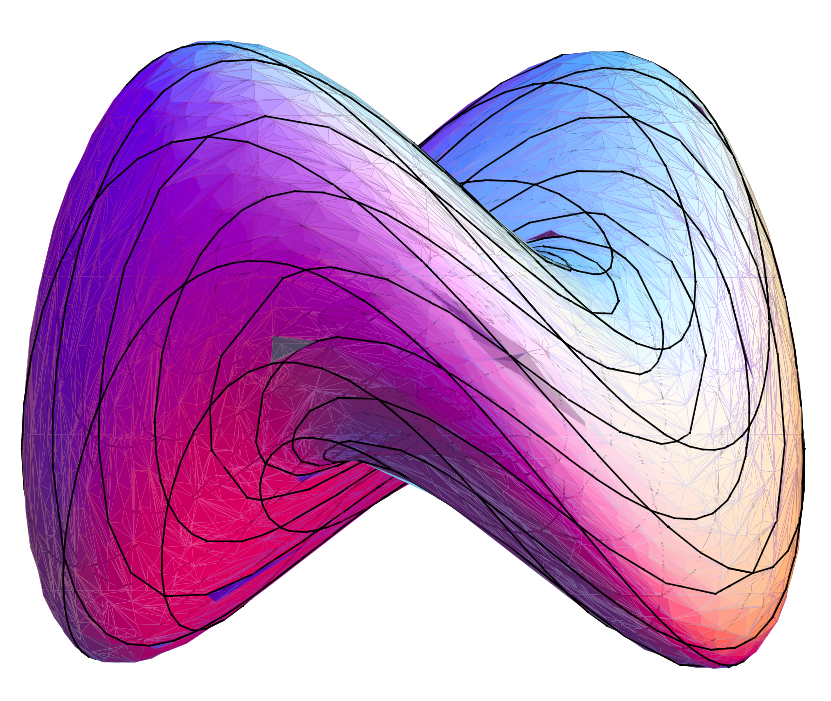}
\caption{Two points of view of the case $\sigma=0,5$ }\label{fig-mixed}
\end{center}
\end{figure}

%
%

\section{Some Remarks on  the heat diffusion}
\label{s-heat}

In this section we briefly discuss the heat diffusion on 3-ARSs.

For a sub-Riemannian manifold, the Laplace operator is defined as the divergence of the horizontal gradient \cite{AgrBarBoscbook,camillo}. The divergence is computed with respect to a given volume, while the horizontal gradient is computed using an orthonormal frame $\{X_1,\ldots,X_m\}$ via the formula grad$_H(\phi)=\sum_i^m X_i(\phi) X_i$.

In particular in the case of 3-ARSs, we have the following.

\bdeff
\label{d-lapl}
Consider a 3-ARS on a smooth manifold $M$. Let $\{X_1,X_2,X_3\}$ be an orthonormal frame defined in an open set $U\subset M$ and let $\mu$ be a smooth volume on $M$.
Then the Laplacian on $U$ is defined as 
$$
\Delta_\mu \phi =\mbox{div$_\mu$(grad}_H(\phi))=\sum_{i=1}^3 \left(X_i^2+\mbox{div$_\mu$}(X_i) X_i\right)\phi
$$
Here div$_\mu$ is the divergence with respect to the volume $\mu$. 
\edeff

\brem
We recall that if  $X=(X^1(x,y,z),X^2(x,y,z),X^3(x,y,z))$ and $\mu=h(x,y,z)dx\,dy\,dz$ then
$$
\mbox{div}_\mu X=\frac{1}{h}\left(\partial_x (h X^1)+\partial_y (h X^2)+\partial_z (h X^3) \right).
$$
\erem

It is easy to check that the definition of $\Delta_\mu$ does not depend on the choice of the orthonormal frame and that  $\Delta_\mu$ is well defined on the whole manifold $M$.

By direct application of the Hormander theorem \cite{hormander} (thanks to the fact that  $\{X_1,X_2,X_3\}$ is bracket generating) and using  a theorem of Strichartz \cite{strichartz}, we have the following.

\bt[Hormander-Strichartz] \label{t-strichartz} Consider a 3-ARS that is complete as metric space. Let $\mu$ be a smooth volume on $M$. Then $\Delta_\mu$ is hypoelliptic and it is essentially self-adjoint on $L^2(M,\mu)$. Moreover the unique solution to the Cauchy problem
\begin{align}
\left\{\ba{l}
(\partial_t-\Delta_\mu)\phi =0\\
\phi(q,0)=\phi_0(q)\in L^2(M,\mu) \cap L^1(M,\mu),\ea\right.
\label{eq-cauchy-gen_H}
\end{align}
on $[0,\infty[\times M$ can be written as 
$$
\phi(q,t)=\int_M\phi_0(\bar q)K_t(q,\bar q) \,\mu(\bar q)
$$
where $K_t(q,\bar q)$ is a positive function defined on $]0,\infty[\times M\times M$ which is smooth, symmetric for the exchange of $q$ and $\bar q$ and such that for every fixed $t,q$, we have  $K_t(q,\cdot) \in L^2(M,\mu) \cap L^1(M,\mu)$.
\et

Theorem \ref{t-strichartz} gives important information on the heat diffusion. Even more, one can relate the heat-kernel asymptotics with the Carnot Caratheodory distance, using the theory developed in {\cite{srneel,BenArous,benarousleandrediag,leandremaj,leandremin}.
For instance a result due to Leandre \cite{leandremaj,leandremin} says that 
\bqn
\label{leandre}
\lim_{t\to0} \Big(-4t \log K_t(q_1,q_2)\Big)=d(q_1,q_2)
\eqn
 In some cases an integral representation for the heat kernel can also be obtained  (see Appendix  \ref{b-heat} for  the case of  nilpotent structures for \Tone\ points,\footnote{
Notice that, in the nilpotent case, if one uses the Lebesgue volume, \Tone\ points are the only interesting ones. Indeed 
the heat kernel for the nilpotent approximation for Riemannian points is well known and  \Ttwo\ points are a particular case of  \Tone\ points.}
 with respect to the Lebesgue volume in $\R^3$).

It should be noticed that the definition of the Laplacian given in Definition \ref{d-lapl} is not completely satisfactory, due to the need of  an external volume $\mu$. One would prefer to define a more intrinsic Laplacian depending only on the 3-ARS. 

An intrinsic choice of volume exists. It is the Riemannian volume $\omega$ associated with the local orthonormal frame $X_1,X_2,X_3$. However 
this volume is well defined only on  $M\setminus\Zz$. See formula \eqref{omega} for its expression using the local representation given by Theorem \ref{t-lokal}.
Hence the Lapacian $\Delta_\omega$ (that we call the intrinsic Laplacian) contains some diverging first order terms and it is well defined only on $M\setminus\Zz$. Using the local representation given by theorem \ref{t-lokal} we obtain,
$$
\Delta_\omega=\partial_x^2+(\al \partial_y+ \beta \partial_z)^2+\nu^2\partial_z^2-\frac{\partial_x (\al\nu )}{\al\nu}\partial_x+ (-\al \frac{\partial_y\nu}{\nu}+\partial_z\beta-\beta \frac{\partial_z (\al\nu )}{\al\nu} )(\al\partial_y+\beta\partial_z)-\frac{\nu^2}{\al}\partial_z\al\,\partial_z.
$$

Theorem \ref{t-strichartz} does not apply to $\Delta_\omega$. This operator is not well defined on the whole manifold.  Theorem \ref{t-strichartz} cannot be applied  even on a connected  component $\Omega$ of $M\setminus\Zz$. Indeed due to the fact that the geodesics can cross the singular set, it happens that in general the 3-ARS restricted to $\Omega$,  is not complete as metric space.

These facts are well known in dimension 2 \cite{camillo}, together with the  fact that $\Zz$ behaves as a barrier for the heat flow. 

In dimension 3 we are going to illustrate that the same phenomenon occurs for the nilpotent structure of \Tone\ points. The fact that $\Zz$ behaves as a barrier for the heat flow is probably true in much more general situations, but this  discussion is out of the purpose of this paper.

\bt 
\label{t-canicola}
On $\R^3$ consider the 3-ARS defined by the following 3 vector fields
$$
\hat X_1(x,y,z)=\left(
\ba{c} 1\\0\\0\ea
\right),~~~
\hat X_2(x,y,z)=\left(
\ba{c} 0\\1\\\cos(\sigma)x\ea
\right),~~~
\hat X_3(x,y,z)=\left(
\ba{c} 0\\0\\\sin(\sigma)x\ea
\right),
$$
where $\sigma\in ]0,\pi/2]$ is a parameter.  The corresponding Riemannian volume (defined on $\R^3\setminus \{x=0\}$) is
$$
\omega=\frac{1}{\sin(\sigma)|x|}dx\,dy\,dz
$$
The intrinsic Laplacian  has the expression 
$$
\Delta_\omega=\partial_x^2 +(\partial_y +\cos(\sigma) x\,\partial_z)^2+\sin(\sigma)^2 x^2 \partial_z^2-\frac{1}{x}\partial_x.
$$
This operator with domain ${\con^\infty_c}(\R^3\setminus\{x=0\})$ is essentially self-adjoint in $L^2(\R^3\setminus\{x=0\})$. Hence it separates in the direct sum of its restrictions to $\R^3\setminus\{x<0\}$ and  $\R^3\setminus\{x>0\}$.

\et

{\bf Proof.} 
\newcommand\Tu{\mathbb{T}}
\newcommand\bna{\begin{eqnarray*}}
\newcommand\ena{\end{eqnarray*}}
Let us make the change of variable in Hilbert space $f=\sqrt{\sin(\sigma)|x|}g$ which is unitary from $L^2(\R^3,\omega)$ to $L^2(\R^3,dx\,dy\,dz)$, so that $(f_1,f_2)_{L^2(\R^3,\omega)}= (g_1,g_2)_{L^2(\R^3,dx\,dy\,dz)}$.

We compute the operator in the new variable:
\bna
\Delta_\omega f&=&\Big(\partial_x^2 +(\partial_y +\cos(\sigma) x\,\partial_z)^2+\sin(\sigma)^2 x^2 \partial_z^2-\frac{1}{x}\partial_x\Big)f=\nn\\
&=&\sqrt{\sin(\sigma)|x|}\Big(\partial_x^2 +(\partial_y +\cos(\sigma) x\,\partial_z)^2+\sin(\sigma)^2 x^2 \partial_z^2-\frac{3}{4 x^2}\Big)g=:\sqrt{\sin(\sigma)|x|}Lg.
\ena
Hence we are left to study the operator $Lg=\Big(\partial_x^2 +(\partial_y +\cos(\sigma) x\,\partial_z)^2+\sin(\sigma)^2 x^2 \partial_z^2-\frac{3}{4 x^2}\Big)g$ on $L^2(\R^3,dx\,dy\,dz)$. By making the Fourier transform in $y$ and $z$, we have $\partial_y\to i\mu$, $\partial_z\to i\nu$, we are left to study the operator
\bna
\hat L^{\mu,\nu}=\partial_x^2 -(\mu +\cos(\sigma) \nu x)^2-\sin(\sigma)^2\nu^2 x^2 -\frac{3}{4 x^2}=\partial_x^2 -V^{\mu,\nu}(x).
\ena
with $V^{\mu,\nu}(x)\geq \frac{3}{4x^2}$. But  in dimension $1$ the operator $-\partial_x^2 +V$ with domain $C^{\infty}_0(]0,+\infty[)$ is essentially self-adjoint on $L^2(]0,+\infty[)$ if $V\geq \frac{3}{4x^2}$ (see \cite{ReedSimon}, Theorem X.10 for the proof of the limit point case at $0$ and Theorem X.8 at $+\infty$). Hence, each operator $L^{\mu,\nu}$ is essentially self-adjoint in $]0,+\infty[$ . 
As a consequence $L^{\mu,\nu}$ is essentially self-adjoint in $]-\infty,+\infty[$ and it separates in the direct sum of its restrictions to    $]-\infty,0[$ and $]0,+\infty[$. By making the inverse Fourier transform, the thesis follows.
\hfill$\Box$

As a direct consequence we have the following
\bc
With the notations of Theorem \ref{t-canicola}, consider the unique solution $\phi$ of the heat equation (according to the self-adjoint extension defined in the previous theorem),
\bqn
\partial_t \phi-\Delta_\omega \phi&=&0\\
\phi(0)&=&\phi_0\in  L^2(\R^3,\omega) \cap L^1(\R^3,\omega)
\eqn
with $\phi_0$ supported in $\R^3\setminus\{x<0\}$. Then, $\phi(t)$ is supported in $\R^3\setminus\{x<0\}$  for any $t\geq 0$.
The same holds for the solution of the Schroedinger equation  or for the solution of the wave equation. 
\ec

Hence formula  \eqref{leandre} does not apply for the diffusion generated by the intrinsic Laplacian. Indeed for \Tone\ points in the nilpotent case, the heat does not flow through $\{x=0\}$, while the Carnot Caratheodory distance is finite for every pair of points.

\appendix
\section{Genericity of (G1),(G2),(G3)}
\label{a-genericity}
In this part, we provide a proof of Proposition~\ref{p-generic}. Before that, we first give some basic results on transversality theory.
\subsection{Thom Transversality Theorem}
Let $M$ and $N$ be smooth manifolds and $k\geq 0$ be an integer. Let $x\in M$, $y\in N$ and $C^\infty(M,N,x,y)$ be the set of smooth maps from $M$ to $N$ which send $x$ to $y$. Let $\varphi$ and $\psi$ be local charts of $M$ and $N$ around $x$ and $y$ respectively.

We use $\mathcal{R}_{k}$ to define the following equivalence relation on $C^\infty(M,N,x,y)$: Two functions $f$ and $g$ are equivalent if the functions $\psi\circ f\circ\varphi^{-1}$ and $\psi\circ g\circ\varphi^{-1}$ have the same partial derivatives at any order less or equal to $k$ at $\varphi(x)$.
\begin{remark}
Notice that ${\mathcal R}_{k}$ is independent of the choice of the charts $\varphi$ and $\psi$. 
\end{remark}
In order to state the Thom Transversality Theorem, we need to list some definitions.
\begin{definition}
Let $M$ and $N$ smooth manifolds. A jet at order $k$ from $M$ to $N$ is a triplet $(x,y,u)$ where $x\in M$, $y\in N$ and $u$ is an equivalence class (for\ ${\mathcal R}_{k}$) of the functions $C^\infty(M,N,x,y)$.

A $k$-th order jet space from $M$ to $N$ denoted by $J^k(M,N)$ is the set of jets at order $k$ from $M$ to $N$.
\end{definition}

\begin{proposition}[\cite{hirsch}]
Let $M$ and $N$ smooth manifolds and $k\geq 0$ an integer. Then $J^k(M,N)$ has a structure of smooth manifold.
\end{proposition}

\begin{definition}
\label{def_general}
Let $M$ and $N$ smooth manifolds. 
\begin{itemize}
\item[$(i)$] We say that a subset of $C^\infty(M,N)$ is residual (and hence dense) if it is an intersection of open dense subsets of  $C^\infty(M,N)$ endowed with the $C^\infty$-Whitney topology.
\item[$(ii)$] We say that $f\in C^\infty(M,N)$ is transverse to a smooth submanifold $S$ of $N$ at $x\in M$ if either $f(x)\notin S$ or $y:=f(x)\in S$ and $Df(x)(T_x M)+T_y S=T_y N$. If $f$ is transverse to $S$ at every point of $M$ then we say that $f$ is transverse to $S$ and we denote it by $f\pitchfork S$. Moreover $f^{-1}(S)$ is a submanifold of $M$ with the same codimension as $S$.
\item[$(iii)$] If $f\in C^\infty(M,N)$ then its $k$-jet extension $J^k f$ is the smooth map from $M$ to $J^k(M,N)$ which assigns to every $x\in M$ the jet of $f$ of order $k$ at $x$.
\end{itemize}
\end{definition}

\begin{theorem}[Thom Transversality Theorem, \cite{hirsch}, Page 82]
\label{Thom}
Let $M, N$ smooth manifolds and $k\geq 1$ an integer. If $S_1,\cdots, S_r$ are smooth submanifolds of $J^k(M,N)$ then the set
\begin{align}
\label{Transv_cond}
\{f\in C^\infty(M,N): J^k f \pitchfork S_i~\text{for}~i=1,2,\cdots,r\},
\end{align}
is residual in the $C^\infty$-Whitney topology.
\end{theorem}

If codim$\ S_i > \dim M$ then $J^k f\pitchfork S_i$ means that $J^k f(M)\cap S_i=\emptyset$. Hence, we have:
 
\begin{corollary}
\label{ba_corollaire}
Assume that $codim\ S_i > \dim M$ for $i=1,\cdots,r$ and $k\geq 1$. Then the set 
\begin{equation}
\{f\in C^\infty(M,N): J^k f\cap S_i=\emptyset~\text{for}~i=1,2,\cdots,r\},
\end{equation}
is residual in the $C^\infty$-Whitney topology.
\end{corollary}

By using Item~$(ii)$ of Definition~\ref{def_general} and Theorem~\ref{Thom} , we have the following:

\begin{corollary}
\label{c_manif}
For every $f$ in the residual set defined in Theorem~\ref{Thom}, the inverse images $\tilde{S_i}:=(J^k f)^{-1}(S_i)$ is a smooth submanifold of $M$ and $codim\ S_i=codim\ \tilde{S}_i$ for $i=1,\cdots,r$.
\end{corollary} 

\subsection{Proof of Proposition~\ref{p-generic}}
Here $M$ is a fixed $3$-dimensional smooth manifold. We use $N$ and $C^\infty(M,N)$ to denote the smooth manifold $\underset{q\in M}\bigcup N_q$ of dimension $12$, 
where $N_q:=(T_q M)^3$, and the set of smooth maps from $M$ to $N$ which assign to every $q\in M$ an element of $N_q$, respectively.
Let recall the conditions $\textbf{(G1)}$, $\textbf{(G2)}$ and $\textbf{(G3)}$:

\hspace{.1cm}{(G1)}  dim($\bD(q))\geq 2$
and $\bD(q)+[\bD,\bD](q)=T_qM$  for every $q\in M$;

\hspace{.1cm}{(G2)} $\Zz$ is an embedded two-dimensional 
submanifold of $M$;

\hspace{.1cm}{(G3)}  the points where $\bD(q)=T_q\Zz$ are isolated. 

\medskip

\noindent\textbf{Proof of Property $\textbf{(G1)}$:} Let us now prove the first part. For this purpose, consider the smooth submanifold of $J^1(M,N)$ of codimension $9$ defined as follows:
$$
S_1:=\{J^1(X_1,X_2,X_3)(q)\in J^1(M,N): (X_1(q), X_2(q), X_3(q))=0\in N_q\}.
$$

Then by using Corollary~\ref{ba_corollaire}, we obtain that
$$
\mathcal{O}_1:=\{(X_1,X_2,X_3)\in C^\infty(M,N): J^1(X_1,X_2,X_3)(M)\cap S_1=\emptyset\},
$$
is a residual subset of $C^\infty(M,N)$ endowed with the $C^\infty$-Whitney topology. We next define the set
$$
S_2:=\left\{J^1(X_1,X_2,X_3)(q)\in J^1(M,N):
\begin{array}{l}
(X_1(q),X_2(q),X_3(q))\neq 0\in N_q,\\ 
~~~~X_1(q)\wedge X_2(q)=0,\\
~~~~X_1(q)\wedge X_3(q)=0,\\
~~~~X_2(q)\wedge X_3(q)=0.
\end{array}
\right\},
$$
and we easily prove that it is a smooth submanifold of $J^1(M,N)$ of codimension strictly greater than $3$.
Therefore we have that the set

$$
\mathcal{O}_2:=\{(X_1,X_2,X_3)\in C^\infty(M,N): J^1(X_1,X_2,X_3)(M)\cap S_2=\emptyset\},
$$
is a residual subset of $C^\infty(M,N)$ in the $C^\infty$-Whitney topology. Thus, for every $(X_1,X_2,X_3)\in\mathcal{O}:=\mathcal{O}_1\cap\mathcal{O}_2$, we have that $J^1(X_1,X_2,X_3)(q)\notin S_1$ and $J^1(X_1,X_2,X_3)(q)\notin S_2$. Hence we conclude the first part of $(\textbf{G1})$. 

We next prove the second step. As above, we define the following subset of $J^1(M,N)$:
$$
S:=\left\{J^1(X_1,X_2,X_3)(q)\in J^1(M,N):
\begin{array}{l}
~~(X_1(q)\wedge X_2(q),X_1(q)\wedge X_3(q),X_2(q)\wedge X_3(q))\neq 0,\\
\det(X_1(q),X_2(q),X_3(q))=0,~\text{and for}~1\leq i<j\leq 3,\\~~~~~~~~
\det(X_i(q),X_j(q),[X_1,X_2](q))=0,\\~~~~~~~~
\det(X_i(q),X_j(q),[X_1,X_3](q))=0,\\~~~~~~~~
\det(X_i(q),X_j(q),[X_2,X_3](q))=0.
\end{array}
\right\}.
$$
By the same strategy as in the first part, we can easily see that $S$ is a smooth submanifold of $J^1(M,N)$ with codimension $4$. Thanks to Corollary~\ref{ba_corollaire}
$$
\mathcal{P}:=\{(X_1,X_2,X_3)\in C^\infty(M,N): J^1(X_1,X_2,X_3)(M)\cap S=\emptyset\},
$$
is a residual subset of $C^\infty(M,N)$ endowed with the $C^\infty$-Whitney topology. Let us denote by $\mathcal A$ the residual set $\mathcal{P}\cap\mathcal{O}$, where $\mathcal O$ is defined in the previous part. Hence we conclude that 
Span$\{X_1(q),X_2(q),X_3(q),[X_1,X_2](q),[X_1,X_3](q),[X_2,X_3](q)\}=T_{q}M$, for~every $(X_1,X_2,X_3)\in\mathcal{A}$ and for every $q\in M$ such that $\det(X_1(q),X_2(q),X_3(q))=0$. This proves the second part of $(\textbf{G1})$.

\medskip

\noindent\textbf{Proof of Property $\textbf{(G2)}$:} For every $(X_1,X_2,X_3)\in C^\infty(M,N)$, we use $\psi_{(X_1,X_2,X_3)}$ and $\bar{S}$ respectively to denote the smooth map
$$
q\in M\rightarrow \det(X_1(q),X_2(q),X_3(q))\in\mathbb{R}.
$$
and the set
$$
\left\{J^1(X_1,X_2,X_3)(q)\in J^1(M,N):
\begin{array}{l}
(X_1(q)\wedge X_2(q),X_1(q)\wedge X_3(q),X_2(q)\wedge X_3(q))\neq 0,\\
~~~~~~~~~~~~~~~~~~~~~D\psi_{(X_1,X_2,X_3)}(q)=0,\\
~~~~~~~~~~~~~~~~~~~~~~~~\psi_{(X_1,X_2,X_3)}(q)=0.
\end{array}
\right\}.
$$
Then $\bar{S}$ is a smooth submanifold of $J^1(M,N)$ of codimension $4$ which implies that the set
$$
\bar{\mathcal{O}}:=\{(X_1,X_2,X_3)\in C^\infty(M,N): J^1(X_1,X_2,X_3)(M)\cap \bar{S}=\emptyset\},
$$
is a residual subset of $C^\infty(M,N)$ in the $C^\infty$-Whitney topology. Let $\bar{\mathcal A}:=\bar{\mathcal O}\cap\mathcal{O}$ and $(X_1,X_2,X_3)\in~\bar{\mathcal{A}}$. Thus, for every $q\in M$ such that $\psi_{(X_1,X_2,X_3)}(q)=0$ we obtain that $D\psi_{(X_1,X_2,X_3)}(q)\neq 0$. This implies that $\psi_{(X_1,X_2,X_3)}$ is transverse to $\{0\}\subset\mathbb{R}$. Hence the inverse image
$$
\{q\in M: \det(X_1(q),X_2(q),X_3(q))=0\},
$$ 
is an embedded submanifold of $M$ of codimension $1$. This proves Property $(\textbf{G2})$.

\medskip

\noindent\textbf{Proof of Property $\textbf{(G3)}$:}  We use the same techniques as previously  by considering the following smooth submanifold of $J^1(M,N)$ of codimension $3$:
$$
\tilde{S}:=\left\{J^1(X_1,X_2,X_3)(q)\in J^1(M,N):
\begin{array}{l}
~~(X_1(q)\wedge X_2(q),X_1(q)\wedge X_3(q),X_2(q)\wedge X_3(q))\neq 0,\\
\det(X_1(q),X_2(q),X_3(q))=0,~\text{and for}~1\leq i\leq 3,\\~~~~~~~~
\det(DX_1(q)X_i(q),X_2(q),X_3(q))\\~~~~~~~~+\det(X_1(q),DX_2(q)X_i(q),X_3(q))\\~~~~~~~~+ \det(X_1(q),X_2(q),DX_3(q)X_i(q))=0.
\end{array}
\right\}.
$$
Then by Theorem~\ref{Thom},
$$
\tilde{\mathcal O}:=\{(X_1,X_2,X_3)\in C^\infty(M,N):J^1(X_1,X_2,X_3)\pitchfork\tilde{S}\}\cap\mathcal{\bar A},
$$
is a residual subset of $C^\infty(M,N)$ in the $C^\infty$-Whitney topology. Now we consider the smooth maps $(X_1,X_2,X_3)\in \tilde{\mathcal O}$. Then, by Corollary~\ref{c_manif}, $J^1(X_1,X_2,X_3)^{-1}(\tilde S)$ is a smooth submanifold of $M$ of codimension ~$3$, i.e., it is formed by isolated points. On the other hand notice that
$$
J^1(X_1,X_2,X_3)^{-1}(\tilde S)=\left\{q\in M:
\begin{array}{l}
\det(X_1(q),X_2(q),X_3(q))=0,~\text{and for}~1\leq i\leq 3,\\~~~~~~~~
\det(DX_1(q)X_i(q),X_2(q),X_3(q))\\~~~~~~~~+\det(X_1(q),DX_2(q)X_i(q),X_3(q))\\~~~~~~~~+ \det(X_1(q),X_2(q),DX_3(q)X_i(q))=0.
\end{array}
\right\},
$$
is the set of points $q\in \mathcal Z$ such that $span\{X_1(q),X_2(q),X_3(q)\}=T_{q}\mathcal{Z}$. Here $\mathcal{Z}$ is the two dimensional embedded submanifold $\{q\in M:\det(X_1(q),X_2(q),X_3(q))=0\}$. Hence $(\textbf{G3})$ is proved.

\section{Explicit expressions of heat kernels}\label{b-heat}
In this section, we consider the nilpotent structures of type-$1$ points, and the Laplacian 
\begin{equation}
\Delta_{{\mathrm L}}:=\partial_x^2+(\partial_ y+x\cos(\sigma)\partial_z)^2+(x\sin(\sigma)\partial_z)^2.
\end{equation}
with respect to the Lebesgue volume $dv=dx\,dy\,dz$ in $\mathbb{R}^3$. In order to give an explicit formula of the associated heat kernel, we first introduce the following intermediate functions:
$$
F(\nu,t):=-t\sin^2\sigma-\frac{\tanh(\nu t)\cos^2\sigma}{\nu},~~~G(\nu,t):=-\cos\sigma\tanh(\nu t),
$$
defined on $(\mathbb{R}\setminus\{0\})\times ]0,+\infty[$. Observe that $F(\nu,t)<0$ for every $\nu\neq 0$ and $t>0$. This comes from the fact that $\frac{\tanh(x)}{x}>0$ for every $x\neq 0$.

Let us also define the next function defined on $]0,+\infty[\times \mathbb{R}^3\times \mathbb{R}^3\times(\mathbb{R}\setminus\{0\})$: 
\begin{eqnarray*}
I(t;x,y,z;\bar x,\bar y,\bar z;\nu):&=&\frac{1}{(2\pi)^2}\cos\left(\nu(z-\bar z)-(x+\bar x)(y-\bar y)\frac{G(\nu,t)}{2 F(\nu,t)}\right)\\
&&\times \exp\left(x \bar x\left(\frac{\nu}{\sinh(2\nu t)}-\frac{G^2(\nu,t)}{2 F(\nu,t)}\right)-(x^2+\bar{x}^2)\left(\frac{\nu}{2\tanh(2\nu t)}+\frac{G^2(\nu,t)}{4 F(\nu,t)}\right)\right)\\
&&\times \exp\left(\frac{(y-\bar y)^2}{4 F(\nu,t)}\right)\times \left(\frac{-\nu}{2 F(\nu,t)\sinh(2\nu t)}\right)^{\frac{1}{2}}.
\end{eqnarray*}


Thus, thanks to Theorem~\ref{t-strichartz} we have the following:

\begin{theorem}
\label{solution-cauchy}
The unique solution of the Cauchy problem
\begin{align}
\left\{\ba{l}
(\partial_t-\Delta_{\mathrm L})\phi =0\\
\phi(x,y,z,0)=\phi_0(x,y,z)\in L^2(\mathbb{R}^3,dv)\cap L^1(\mathbb{R}^3,d v),\ea\right.
\label{eq-cauchy-L}
\end{align} 
defined on $\mathbb{R}^3\times [0,\infty[$ is of the form
$$
\phi(x,y,z,t)=\int_{\mathbb{R}^3}\phi_0(\bar x,\bar y,\bar z)K_t(x,y,z;\bar x,\bar y,\bar z) \,
d{\bar x}\ d{\bar y} \ d{\bar z},
$$
where
\begin{align}
\label{NChaleur}
K_t(x,y,z;\bar x,\bar y,\bar z)=\int_{\mathbb{R}}{I(t,x,y,z;\bar x,\bar y,\bar z,\nu)}\ d\nu.
\end{align}
\end{theorem}
{\bf Proof.}
Let $\phi$ and $\hat{\phi}$  the solution of Problem~(\ref{eq-cauchy-L}) and its Fourier transform. Applying the inverse Fourier transform only on $y$ and $z$, we get that
$$
\phi(x,y,z,t)=\frac{1}{(2\pi)^2}\int_{\mathbb{R}^2} \exp(i(\mu y+\nu z))\hat{\phi}(x,\mu,\nu,t)\ d\mu \ d\nu.
$$
Thus, we easily prove that Problem~(\ref{eq-cauchy-L}) is equivalent to the following:
\begin{align}
\left\{\ba{l}
\partial_t{\hat{\phi}(x,\mu,\nu,t)}=(\partial_x^2-(\nu x+\mu\cos\sigma)^2-\mu^2\sin^{2}\sigma)\hat{\phi}(x,\mu,\nu,t),\\
\hat{\phi}(x,\mu,\nu,0)={\hat{\phi}}_{0}(x,\mu,\nu)\in L^2(\mathbb{R}^3,\mathbb{R})\cap L^1(\mathbb{R}^3,\mathbb{R}).\ea\right.
\label{eq-cauchy-t-fourier}
\end{align}
Hence, by making the change of variable $x\rightarrow\gamma=x+\frac{\mu}{\nu}\cos\sigma$ ($\nu\neq 0$), Problem~(\ref{eq-cauchy-t-fourier}) becomes:
\begin{align}
\left\{\ba{l}
\partial_{t}\ {\bar{\phi}^{\mu,\nu}(\gamma,t)}=\left({\partial_\gamma}^{2}-\nu^2\gamma^2-\mu^2\sin^2\sigma\right)\bar{\phi}^{\mu,\nu}(\gamma,t),\\
\bar{\phi}^{\mu,\nu}(\gamma,0)=\bar{\phi}_{0}^{\mu,\nu}(\gamma),
\ea\right.
\label{eq-cauchy-t-fourier-var-change}
\end{align}
where
$$
\bar{\phi}^{\mu,\nu}(\gamma,t):=\hat{\phi}(\gamma-\frac{\mu}{\nu}\cos\sigma,\mu,\nu,t)~~\text{and}~~\bar{\phi}_{0}^{\mu,\nu}(\gamma):={\hat{\phi}}_{0}(\gamma-\frac{\mu}{\nu}\cos\sigma,\mu,\nu).
$$
In the sequel, we use $\psi(\gamma,t)$ to denote the solution of Problem~(\ref{eq-cauchy-t-fourier-var-change}). First remark that the eigenvalues and the associated eigenfunctions of the operator $\partial_{\gamma}^{2}-\nu^2\gamma^2-\mu^2\sin^2\sigma$\ ($\nu\neq 0$ ) on $\mathbb{R}$ are respectively given by 
$$
E_n=-2\nu(n+\frac{1}{2})-\mu^2 \sin^2\sigma,
$$
and
$$
\phi_n^{\nu}(\gamma):=\frac{1}{\sqrt{2^n \fact{n}}}\left(\frac{\nu}{\pi}\right)^{\frac{1}{4}} \exp(-\frac{\nu\gamma^2}{2})\ H_n(\gamma\sqrt{\nu})~~\text{where}~~ H_n(\gamma):=(-1)^{n} e^{\gamma^2}\frac{d^n}{d\gamma^n} \exp(-\gamma^2),~~n=0,1,\cdots.
$$
Since the sequence $\{\phi_n^{\nu}\}_{n}$ is an orthonormal basis of $L^2(\mathbb{R})$ then there exists a sequence of functions $\{C_n(\cdot)\}_{n}$ such that $\psi(\gamma,t)=\sum_{n\geq 0}{C_n(t)\phi_n^{\nu}(\gamma)}$. According to Problem~(\ref{eq-cauchy-t-fourier-var-change}), we easily obtain that $\dot{C_n}(t)=E_n C_n(t)$ which implies that $C_n(t)=\exp(t E_n)C_n(0)$ for $n\in\mathbb{N}$ and $t\geq 0$. The fact that 
$$
C_{n}(0)=\int_{\mathbb{R}} \psi_{0}(\bar\gamma)\phi_n^{\nu}(\bar\gamma) d{\bar\gamma},
$$
implies after simple computations that
$$
\psi(\gamma,t)=\int_{\mathbb{R}} \psi_{0}(\bar\gamma){Q_{t}}^{\mu,\nu}(\gamma,\bar\gamma)\ d{\bar\gamma}~~\text{with}~~{Q_{t}}^{\mu,\nu}(\gamma,\bar\gamma)=\left(\sum_{n\geq0}{\exp(t E_n)\ \phi_n^{\nu}(\gamma)\ \phi_n^{\nu}(\bar\gamma)}\right).
$$
Thus, we obtain that
$$
{Q_{t}}^{\mu,\nu}(\gamma,\bar\gamma)=\left(\frac{\nu}{\pi}\right)^{\frac{1}{2}}\exp\left(-\frac{\nu}{2}(\gamma^2+\bar\gamma^2+2t)-t\mu^2\sin^2\sigma \right)\ \sum_{n\geq0}{\exp(-2t\nu n)}\ \frac{1}{2^n\fact{n}}\ H_n(\gamma\sqrt{\nu})\ H_n(\bar\gamma\sqrt{\nu}).
$$
Let us denote $w=\exp(-2 t\nu)$. By the Mehler's formula we get that
$$
\sum_{n\geq0}{\exp(-2t\nu n)}\ \frac{1}{2^n\fact{n}}\ H_n(\gamma\sqrt{\nu})\ H_n(\bar\gamma\sqrt{\nu})=(1-w^2)^{-\frac{1}{2}}\ \exp\left(\frac{2\nu \gamma\bar{\gamma} w-\nu(\gamma^2+\bar\gamma^2)w^2}{1-w^2}\right),
$$
which implies that
$$
{Q_{t}}^{\mu,\nu}(\gamma,\bar\gamma)=\left(\frac{\nu}{\pi}\right)^{\frac{1}{2}}\exp\left(-\frac{\nu}{2}(\gamma^2+\bar\gamma^2+2t)-t\mu^2\sin^2\sigma \right)\ (1-w^2)^{-\frac{1}{2}}\ \exp\left(\frac{2\nu \gamma\bar{\gamma} w-\nu(\gamma^2+\bar\gamma^2)w^2}{1-w^2}\right) .
$$
After some algebraic computations, we deduce that 
$$
{Q_{t}}^{\mu,\nu}(\gamma,\bar\gamma)=\exp(-t\mu^2\sin^2\sigma)\ \left(\frac{\nu w}{\pi (1-w^2)}\right)^{\frac{1}{2}}\ \exp\left(\frac{\nu\gamma\bar\gamma(w-1)}{w+1}-\frac{\nu(\gamma-\bar\gamma)^2(w^2+1)}{2(1-w^2)}\right).
$$
Let us remark that
$$
\frac{w}{1-w^2}=\frac{1}{2\sinh(2\nu t)},~~ \frac{w-1}{1+w}=\frac{1}{\sinh(2\nu t)}-\frac{1}{\tanh(2\nu t)}=-\tanh(\nu t)~~\text{and}~~\frac{w^2+1}{1-w^2}=\frac{1}{\tanh(2\nu t)}.
$$
Hence, we conclude that
$$
{Q_{t}}^{\mu,\nu}(\gamma,\bar\gamma)=\left(\frac{\nu}{2\pi\sinh(2\nu t)}\right)^{\frac{1}{2}}\  \exp\left(-\left(t\mu^2\sin^2 \sigma+\frac{\nu(\gamma-\bar\gamma)^2}{2\tanh(2\nu t)}+\nu\tanh(\nu t)\gamma\bar\gamma\right)\right).
$$
Since $\psi_{0}(\bar\gamma)=\hat{\phi}_{0}(\bar\gamma-\frac{\mu}{\nu}\cos\sigma,\mu,\nu)$, then we have that
$$
\psi_{0}(\bar\gamma)=\int_{\mathbb{R}^2}{\exp(-i\mu\bar y)\ \exp(-i\nu \bar z)\ \phi_0(\bar\gamma-\frac{\mu}{\nu}\cos\sigma,\bar y,\bar z)}\ d\bar y\ d\bar z,
$$
which implies that
\begin{equation}
\label{F1}
\bar{\phi}^{\mu,\nu}(\gamma,t)=\int_{\mathbb{R}}\left(\int_{\mathbb{R}^2}{\exp(-i\mu\bar y)\ \exp(-i\nu \bar z)\ \phi_0(\bar\gamma-\frac{\mu}{\nu}\cos\sigma,\bar y,\bar z)}\ d\bar y\ d\bar z\right){{Q_{t}}^{\mu,\nu}(\gamma,\bar\gamma)}\ d\bar\gamma.
\end{equation}
By the fact that $\bar{\phi}^{\mu,\nu}(\gamma,t)=\hat{\phi}(\gamma-\frac{\mu}{\nu}\cos\sigma,\mu,\nu,t)$ and by making the inverse Fourier transform we deduce that
\begin{equation}
\label{F2}
\phi(x,y,z,t)=\frac{1}{(2\pi)^2}\int_{\mathbb{R}^2}{\exp(i\mu y)\ \exp(i\nu z)\ \bar{\phi}^{\mu,\nu}\left(x+\frac{\mu \cos\sigma}{\nu},t\right)}\ d\mu\ d\nu.
\end{equation}
Hence by using the change of variable $\bar\gamma\rightarrow \bar x=\bar\gamma-\frac{\mu \cos\sigma}{\nu}$, and replacing Eq.~(\ref{F1}) with $\gamma=x+\frac{\mu \cos\sigma}{\nu}$ in Eq.~(\ref{F2}), we easily conclude that
$$
\phi(x,y,z,t)=\int_{\mathbb{R}^3} K_t(x,y,z;\bar x,\bar y,\bar z)\ \phi_0(\bar x,\bar y,\bar z) \,
d{\bar x}\ d{\bar y} \ d{\bar z},
$$
where
$$
K_t(x,y,z;\bar x,\bar y,\bar z)=\frac{1}{(2\pi)^2}\int_{\mathbb{R}^2}{Q_{t}}^{\mu,\nu}(x+\frac{\mu}{\nu}\cos\sigma,\bar{x}+\frac{\mu}{\nu}\cos\sigma)\ \exp(i\mu(y-\bar y))\ \exp(i\nu(z-\bar z))\,
d{\mu}\ d{\nu}.
$$
Since


%
\begin{eqnarray*}
{Q_{t}}^{\mu,\nu}(x+\frac{\mu}{\nu}\cos\sigma,\bar{x}+\frac{\mu}{\nu}\cos\sigma)&=& \exp(-t\mu^2\sin^2 \sigma)\ \left(\frac{\nu}{2\pi\sinh(2\nu t)}\right)^{\frac{1}{2}}\\
&&\times \exp\left(-\left(\frac{\nu(x-\bar x)^2}{2\tanh(2\nu t)}+\nu\tanh(\nu t)(x+\frac{\mu}{\nu}\cos\sigma)(\bar x+\frac{\mu}{\nu}\cos\sigma)\right)\right),
\end{eqnarray*}
then by defining
$$
A(\nu,t):=\int_{\mathbb{R}}\exp(i\mu(y-\bar y))\exp\left(-t\mu^2\sin^2\sigma-\nu\tanh(\nu t)\left(\frac{\mu(x+\bar x)\cos\sigma}{\nu}+\frac{\mu^2 \cos^2 \sigma}{\nu^2}\right)\right)\ d\mu,
$$
we deduce that
\begin{align}
\label{heat_kernel}
K_t(x,y,z;\bar x,\bar y,\bar z)=\frac{1}{(2\pi)^2} \int_{\mathbb{R}} \exp(i\nu(z-\bar z))\left(\frac{\nu}{2\pi\sinh(2\nu t)}\right)^{\frac{1}{2}} \exp\left(-\frac{\nu(x-\bar x)^2}{2\tanh(2\nu t)}-\nu\tanh(\nu t)x\bar x\right) A(\nu,t)\ d\nu.
\end{align}
Let us now give a better formula of $A(\nu,t)$. By straightforward computations, we observe that
$$
A(\nu,t)=\int_{\mathbb R} \exp(i\mu(y-\bar y))\exp\left(F(\nu,t)\mu^2+(x+\bar x)G(\nu,t)\mu\right)\ d\mu.
$$
Therefore by making the change of variable $\mu\rightarrow \xi=\mu+\frac{(x+\bar x)\ G(\nu,t)}{2 F(\nu,t)}$, we easily have that
$$
A(\nu,t)=\exp\left(-\frac{(x+\bar x)^2\  G^2(\nu,t)}{4\ F(\nu,t)}\right)\int_{\mathbb R} \exp\left(i\ (y-\bar y)\left(\xi-\frac{(x+\bar x)\ G(\nu,t)}{2\ F(\nu,t)}\right)\right)\ \exp\left(F(\nu,t)\xi^2\right)\ d\xi.
$$
Since
$$
\int_{\mathbb R} \exp(i\xi(y-\bar y))\exp\left(F(\nu,t)\xi^2\right)\ d\xi=\left(\frac{-\pi}{F(\nu,t)}\right)^{\frac{1}{2}}\ \exp\left(\frac{(y-\bar y)^2}{4 F(\nu,t)}\right),
$$
we conclude that
\begin{align}
\label{Exact_value}
A(\nu,t)=\left(\frac{-\pi}{F(\nu,t)}\right)^{\frac{1}{2}}\ \exp\left(\frac{(y-\bar y)^2}{4 F(\nu,t)}-\frac{(x+\bar x)^2 \ G^2(\nu,t)}{4\ F(\nu,t)}-i\ \frac{(x+\bar x)\ (y-\bar y)\ G(\nu,t)}{2\ F(\nu,t)}\right).
\end{align}
Hence, we obtain Eq.$(\ref{NChaleur})$ of Theorem~\ref{solution-cauchy} by replacing Eq.(\ref{Exact_value}) in Eq.(\ref{heat_kernel}).
$\Box$\hfill

Observe that the Baouendi-Goulaouic operator which is defined in~(\ref{Baouendi_Goulaouic}) corresponds to $\Delta_{\mathrm{L}}$ in the case where $\sigma=\frac{\pi}{2}$. Hence, according to the previous theorem we have the following
\begin{corollary}
The heat kernel associated with the Baouendi-Goulaouic operator is given by
$$
K_t(x,y,z;\bar x,\bar y,\bar z):=\displaystyle{\int_{\mathbb{R}}}{I_{B}(t;x,y,z;\bar x,\bar y,\bar z;\nu)}\ d\nu,
$$
where 
\begin{eqnarray*}
I_{B}(t;x,y,z;\bar x,\bar y,\bar z;\nu):&=&\frac{1}{(2\pi)^2}\ cos(\nu(z-\bar z))\times \left(\frac{\nu}{2t\sinh(2\nu t)}\right)^{\frac{1}{2}}\\&&
\times \exp\left(\frac{\nu\ x\bar{x}}{\sinh(2\nu t)}-\frac{\nu\ (x^2+\bar x^2)}{2 \tanh(2\nu t)}-\frac{(y-\bar y)^2}{4 t}\right).
\end{eqnarray*}
\end{corollary}

Let us now consider the case where $\sigma=0$. Then we obtain the well-known Heisenberg-operator
$\partial_x^2+(\partial_ y+x\partial_z)^2
$
in $\mathbb{R}^3$. Hence, we get the next corollary.
\begin{corollary}
The heat kernel associated with the Heisenberg-operator is given by
$$
K_t(x,y,z;\bar x,\bar y,\bar z):=\displaystyle{\int_{\mathbb{R}}}{I_{H}(t;x,y,z;\bar x,\bar y,\bar z;\nu)}\ d\nu,
$$
where 
\begin{eqnarray*}
I_{H}(t;x,y,z;\bar x,\bar y,\bar z;\nu):&=&\frac{1}{(2\pi)^2}\ cos\left(\nu\left((z-\bar z)-\frac{(x+\bar x)(y-\bar y)}{
2}\right)\right)\times \left(\frac{\nu}{2 \sinh(\nu t)}\right)\\&&
\times \exp\left(-\frac{\nu}{4 \tanh(\nu t)}\left((x-\bar x)^2+(y-\bar y)^2\right)\right).
\end{eqnarray*}
\end{corollary}

\bigskip
{\bf Acknowledgement }
We are deeply thankful to Andrei for many illuminating discussions. We are thankful to Bronislaw Jakubczyk for some lecture notes on transversality theory that he wrote and gave to us, on which appendix A is based.

\bibliography{biblio-3D}

\bibliographystyle{plain}

\end{document}